\documentclass{article}
\usepackage{amssymb}

\usepackage{amsmath}
\usepackage{float}

\begin{document}

\title{Integrals containing the logarithm of the Airy Function $Ai^{\,\prime
}(x)$}
\author{Bernard J. Laurenzi \\
Department of Chemistry\\
The State Univrsity of New York at Albany}
\date{December 12, 2017 }
\maketitle

\begin{abstract}
Integrals occurring in Thomas-Fermi theory which contains the logarithm of
the Airy function $Ai^{\,\prime }(x)$ have been obtained in terms of
analytical expressions.
\end{abstract}

\section{An integral containing $Ai^{\,\prime }(x)\,ln(Ai^{\,\prime }(x))$}

\bigskip

Integrals\thinspace of the kind%
\begin{equation*}
\mathfrak{I}^{\,(\alpha )}\,\mathfrak{=}\int_{0}^{\infty }\left( \frac{%
Ai^{\,\prime }(x)}{Ai^{\,\prime }(0)}\right) ^{\alpha }\ln \left( \frac{%
Ai^{\,\prime }(x)}{Ai^{\,\prime }(0)}\right) dx,
\end{equation*}%
which contain the logarithm of the Airy function $Ai^{\,\prime }(x)$ \cite%
{airynbs} and which occur in Thomas-Fermi theory, a semiclassical quantum
mechanical theory for the electronic structure of many-body systems, \cite%
{mathphy} have long resisted expression in closed form or accurate values in
terms of classical mathematical constants or analytic special functions. \
Approximate analytic expressions for $\mathfrak{I}^{\,(1)}$ accurate to
f\thinspace ive decimal places have been given \cite{laur}. \ Here we give
analytic expressions for $\mathfrak{I}^{\,(\alpha )}$ for orders $\alpha
=1,2 $ in the form of inf\thinspace inite series.

We begin by noting that integration by parts of $\mathfrak{I}^{\,(1)}$
yields the expression%
\begin{equation*}
\mathfrak{\mathfrak{I}}^{\,(1)}\mathfrak{=-}\frac{1}{Ai^{\,\prime }(0)}%
\int_{0}^{\infty }Ai(x)\frac{d}{dx}\left\{ \ln \left( \frac{Ai^{\,\prime }(x)%
}{Ai^{\,\prime }(0)}\right) \right\} dx.
\end{equation*}%
Progress can be made in the calculation of $\mathfrak{\mathfrak{I}}^{\,(1)}$%
\ using the Weierstrass inf\thinspace inite product representation \cite%
{vallee} for the quantity $Ai^{\,\prime }(x)/Ai^{\,\prime }(0)$ i.e. 
\begin{equation}
\frac{Ai^{\,\prime }(x)}{Ai^{\,\prime }(0)}=\prod_{n=1}^{\infty }\left\{ 1+%
\frac{x}{\left| a_{n}^{\prime }\right| }\right\} \exp \left\{ -\frac{x}{%
\left| a_{n}^{\prime }\right| }\right\} ,  \label{eq1}
\end{equation}%
where $a_{n}^{\,\prime }$ are the roots \cite{nbs} of the Airy function $%
Ai^{\,\prime }(x).$ \ The f\thinspace irst few of those roots are contained
in the Table I below.\pagebreak

\begin{center}
\begin{table}[tbp] \centering%
\caption{Roots of  Ai'(x)\label{one}}%
\end{table}%

\begin{tabular}{|c|c|}
\hline
$n$ & $a_{n}^{\prime }$ \\ \hline
$1$ & $-1.0187929716$ \\ \hline
$2$ & $-3.2481975822$ \\ \hline
$3$ & $-4.8200992112$ \\ \hline
$4$ & $-6.1633073556$ \\ \hline
$5$ & $-7.3721772550$ \\ \hline
$6$ & $-8.4884867340$ \\ \hline
$7$ & $-9.5354490524$ \\ \hline
$8$ & $-10.5276603970$ \\ \hline
$9$ & $-11.4750566335$ \\ \hline
$10$ & $-12.3847883718$ \\ \hline
\end{tabular}%
\bigskip
\end{center}

Values of $a_{n}^{\prime }$ for large $n$ \cite{Olver} are given by 
\begin{equation*}
a_{n}^{\prime }=-t^{\,2/3}\{1-\tfrac{7}{48}t^{-2}+\tfrac{35}{288}t^{-4}-%
\tfrac{181228}{207360}t^{-6}+\cdots \},
\end{equation*}%
where 
\begin{equation*}
t=\frac{3}{8}\pi (4n-3).
\end{equation*}%
Using the expressions above we get for the integral in (1) 
\begin{equation*}
\mathfrak{I}^{\,(1)}\mathfrak{=}\frac{1}{Ai^{\,\prime }(0)}%
\sum_{n=1}^{\infty }\frac{1}{\left| a_{n}^{\prime }\right| }\int_{0}^{\infty
}\frac{x\,Ai(x)}{(x+\left| a_{n}^{\prime }\right| )}\,dx.
\end{equation*}%
Employing the dif\thinspace ferential equation which def\thinspace ines the
Airy function i.e. 
\begin{equation}
Ai^{\,\prime \prime }(x)-x\,Ai(x)=0,  \label{eq2}
\end{equation}%
$\mathfrak{I}^{(1)}$ becomes 
\begin{equation*}
\mathfrak{I}^{(1)}\mathfrak{=}\frac{1}{Ai^{\,\prime }(0)}\sum_{n=1}^{\infty }%
\frac{1}{\left| a_{n}^{\prime }\right| }\int_{0}^{\infty }\frac{Ai^{\,\prime
\prime }(x)}{(x+\left| a_{n}^{\prime }\right| )}\,dx.
\end{equation*}%
A two-fold integration by parts of the latter expression together with the
use of the f\thinspace irst two Airy zeta sums $\mathcal{Z}_{k}$ i.e. (cf.
Appendix I) 
\begin{eqnarray*}
\mathcal{Z}_{2} &=&\sum_{n=1}^{\infty }\frac{1}{\left| a_{n}^{\,\prime
}\right| ^{2}}=-Ai(0)/Ai^{\,\,\prime }(0), \\
\mathcal{Z}_{3} &=&\sum_{n=1}^{\infty }\frac{1}{\left| a_{n}^{\,\prime
}\right| ^{3}}=1,
\end{eqnarray*}%
yields 
\begin{equation}
\mathfrak{I^{(1)}=}\frac{2}{Ai^{\,\prime }(0)}\sum_{n=1}^{\infty }\frac{1}{%
\left| a_{n}^{\prime }\right| }\int_{0}^{\infty }\frac{Ai(x)}{(x+\left|
a_{n}^{\prime }\right| )^{3}}\,\,dx.  \label{eq3}
\end{equation}%
To complete this calculation it is essential to obtain accurate values for
the integrals occurring in (3) especially in the case of integrals
containing the higher roots. \ 

\subsection{The Generalized Stieltjes transforms of $Ai(x)$}

The integrals $\mathbf{I}_{k}(a),$ $(a>0)$%
\begin{equation}
\mathbf{I}_{k}(a)=\int_{0}^{\infty }\frac{\,Ai(x)}{(x+a)^{k}}\,dx,
\label{eq4}
\end{equation}%
which can be identif\thinspace ied as the generalized Stieltjes transforms
of the Airy function $Ai(x)$, have magnitudes which decrease rapidly for
large $k$, $\mathbf{I}_{k}(a)$ varying roughly as $\frac{1}{3a^{\,k}}$ (cf.
Appendix II). \ 

\bigskip

As will be seen below, the $\mathbf{I}_{k}(a)$ integrals are interrelated by
recurrence relations. \ If we rewrite $\mathbf{I}_{k}(a)$ as%
\begin{equation*}
\mathbf{I}_{k}(a)=a\int_{0}^{\infty }\frac{\,Ai(x)}{(x+a)^{k+1}}%
\,dx+\int_{0}^{\infty }\frac{\,xAi(x)}{(x+a)^{k+1}}\,dx,
\end{equation*}%
we get%
\begin{equation*}
\mathbf{I}_{k}(a)=a\,\mathbf{I}_{k+1}(a)+\int_{0}^{\infty }\frac{%
\,Ai^{\,\prime \prime }(x)}{(x+a)^{k+1}}\,dx.
\end{equation*}%
Integrating by parts twice, the integral above produces the recurrence
relation%
\begin{equation}
\mathbf{I}_{k}(a)-a\,\mathbf{I}_{k+1}(a)-(k+1)(k+2)\mathbf{I}%
_{k+3}(a)=-Ai^{\prime }(0)/a^{k+1}-(k+1)\,Ai(0)/a^{k+2}.  \label{eq5}
\end{equation}

Using (5) together with the dif\thinspace ferential relations%
\begin{eqnarray*}
\frac{d^{\,\,}\mathbf{I}_{k+1}(a)}{d\,a} &=&-(k+1)\mathbf{I}_{k+2}(a), \\
\frac{d^{\,\,2}\mathbf{I}_{k+1}(a)}{d\,a^{\,2}} &=&(k+1)(k+2)\mathbf{I}%
_{k+3}(a),
\end{eqnarray*}%
the recurrence relation (5) is transformed into the second-order
dif\thinspace ferential equation%
\begin{equation}
\frac{d^{\,\,2}\mathbf{I}_{k+1}(a)}{d\,a^{\,2}}+a\,\mathbf{I}_{k+1}(a)=%
\mathbf{I}_{k}(a)+\frac{Ai^{\,\prime }(0)}{a^{k+1}}+\frac{(k+1)\,Ai(0)}{%
a^{k+2}}.  \label{eq6}
\end{equation}%
In the case where $k=0$ 
\begin{equation*}
\mathbf{I}_{0}(a)=\frac{1}{3}.
\end{equation*}%
The integral $\mathbf{I}_{1}(a)$ which will be shown to be important in the
sequel satisf\thinspace ies the dif\thinspace ferential equation 
\begin{equation}
\frac{d^{\,2}\,\mathbf{I}_{1}(a)}{d\,a^{2}}+a\,\,\mathbf{I}_{1}(a)=\frac{1}{3%
}+\frac{Ai^{\,\prime }(0)}{a}+\frac{Ai(0)}{a^{2}},  \label{eq7}
\end{equation}%
with initial conditions%
\begin{equation*}
\mathbf{I}_{1}(a)|_{a=a_{0}}=\mathbf{I}_{1}(a_{0}),\hspace{0.25in}[\frac{d\,%
\mathbf{I}_{1}(a)}{d\,a}]_{a=a_{0}}=-\,\,\mathbf{I}_{2}(a_{0}),
\end{equation*}%
where $a_{0}$ is any positive number. \ The general solution to (7) is then
given by 
\begin{eqnarray*}
\mathbf{I}_{1}(a) &=&\pi \,\mathbf{I}_{1}(a_{0})\left[ Ai(-a)\,Bi^{\,\prime
}(-a_{0})-Bi(-a)\,Ai^{\,\prime }(-a_{0})\right] \\
&&-\pi \,\mathbf{I}_{2}(a_{0})\left[ Ai(-a)\,Bi(-a_{0})-Bi(-a)\,Ai(-a_{0})%
\right] \\
&&+\frac{\pi }{3}\{Ai(-a)\int_{a_{0}}^{a}Bi(-z)dz-Bi(-a)%
\int_{a_{0}}^{a}Ai(-z)dz\} \\
&&+Ai(-a)\int_{a_{0}}^{a}[-\frac{3^{1/6}}{2}\frac{\Gamma (2/3)}{z}+\frac{\pi 
}{3^{2/3}\Gamma (2/3)\,z^{2}}\,]Bi(-z)\,dz\, \\
&&-Bi(-a)\int_{a_{0}}^{a}[-\frac{3^{1/6}}{2}\frac{\Gamma (2/3)}{z}+\frac{\pi 
}{3^{2/3}\Gamma (2/3)\,z^{2}}\,]Ai(-z)\,dz.
\end{eqnarray*}%
The integrals appearing above in the expression for $\mathbf{I}_{1}(a)$ in
the f\thinspace irst and second instances are given by Mathematica as%
\begin{eqnarray*}
\int Bi(-z)\,dz &=&\tfrac{z}{3^{1/6}\Gamma (2/3)}\,_{1}F_{2}\tbinom{%
1/3\;;\;-\;z^{3}/9}{2/3,\;4/3\;}-\tfrac{3^{2/3}\Gamma (2/3)\left( \frac{z}{2}%
\right) ^{2}}{\pi }\,\,_{1}F_{2}\tbinom{2/3\;;\;-\,\,z{}^{3}/9}{4/3,\;5/3\;},
\\
&& \\
\int Ai(-z)\,dz &=&\tfrac{z}{3^{2/3}\Gamma (2/3)}\,_{1}F_{2}\tbinom{%
1/3\;;\;-\;z^{3}/9}{2/3,\;4/3\;}+\tfrac{3^{1/6}\Gamma (2/3)\left( \frac{z}{2}%
\right) ^{2}\,}{\pi }\,_{1}F_{2}\tbinom{2/3\;;\;-\,z{}^{3}/9}{4/3,\;5/3\;},
\\
&&
\end{eqnarray*}%
where $_{a\,}F_{\,b}$ are the generalized hypergeometric functions. \
Gathering terms we have%
\begin{eqnarray*}
&&\frac{\pi }{3}\left\{
Ai(-a)\int_{a_{0}}^{a}Bi(-z)\,dz-Bi(-a)\int_{a_{0}}^{a}Ai(-z)\,dz\right\} \\
&& \\
&=&\frac{\Gamma (1/3)}{3^{7/6}}\mathfrak{J}_{-}(a)\left[ \left( \frac{a}{2}%
\right) \,\,_{1}F_{2}\tbinom{1/3\;;\;-z{}^{3}/9}{2/3,\;4/3\;}\right]
_{a_{0}}^{a} \\
&&-\frac{\Gamma (2/3)}{3^{5/6}}\mathfrak{J}_{+}(a)\left[ \left( \frac{a}{2}%
\right) ^{2}\,\,_{1}F_{2}\tbinom{2/3\;;\;-z{}^{3}/9}{4/3,\;5/3\;}\right]
_{a_{0}}^{a},
\end{eqnarray*}%
where we have def\thinspace ined the functions $\mathfrak{J}_{\pm }(a)$ and $%
\mathfrak{J}_{\pm }^{\,\prime }(a)$\ as 
\begin{equation*}
\mathfrak{J}_{\pm }(a)=\sqrt{3}Ai(-a)\pm Bi(-a),\hspace{0.25in}\ \mathfrak{J}%
_{\pm }^{\,\prime }(a)=\sqrt{3}Ai^{\,\prime }(-a)\pm Bi^{\,\prime }(-a).
\end{equation*}%
In the second instances, the integrals 
\begin{eqnarray*}
V_{1}(z) &=&\int [-\frac{3^{1/6}}{2}\frac{\Gamma (2/3)}{z}+\frac{\pi }{%
3^{2/3}\Gamma (2/3)\,z^{2}}]\,Bi(-z)\,dz, \\
V_{2}(z) &=&\int [-\frac{3^{1/6}}{2}\frac{\Gamma (2/3)}{z}+\frac{\pi }{%
3^{2/3}\Gamma (2/3)\,z^{2}}]\,Ai(-z)\,dz,
\end{eqnarray*}%
are given by 
\begin{eqnarray*}
V_{1}(z) &=&-\ln (z)-\frac{\pi }{3^{5/6}\Gamma (2/3)^{2}\,z}\,_{1}F_{2}%
\tbinom{-1/3\;;\;-z^{3}/9}{2/3,\;2/3\;} \\
&&+\frac{3^{5/6}\Gamma (2/3)^{2}}{4\pi }\,z\,_{1}F_{2}\tbinom{%
1/3\;;\;-z^{3}/9}{4/3,\;4/3\;} \\
&&+\frac{1}{72}\,z^{3}\,[2\,\cdot \,_{2}F_{3}\tbinom{1,1\;;\;-z^{3}/9}{%
2,\;2,\,\,5/3\;}+\,_{2}F_{3}\tbinom{1,1\;;\;-z^{3}/9}{2,\;2,\,\,7/3\;}], \\
&& \\
V_{2}(z) &=&-\frac{\pi }{3^{4/3}\Gamma (2/3)^{2}\,z}\,_{1}F_{2}\tbinom{%
-1/3\;;\;-z^{3}/9}{2/3,\;2/3\;} \\
&&-\frac{3^{1/3}\Gamma (2/3)^{2}}{4\pi }\,z\,_{1}F_{2}\tbinom{%
1/3\;;\;-z^{3}/9}{4/3,\;4/3\;} \\
&&+\frac{\sqrt{3}}{216}\,z^{3}\,[2\cdot \,\,_{2}F_{3}\tbinom{1,1\;;\;-z^{3}/9%
}{2,\;2,\,\,5/3\;}-\,_{2}F_{3}\tbinom{1,1\;;\;-z^{3}/9}{2,\;2,\,\,7/3\;}].
\end{eqnarray*}%
\ Collecting all terms containing hypergeometric functions we have for $%
\mathbf{I}_{1}(a)$ 
\begin{eqnarray*}
\mathbf{I}_{1}(a) &=&\frac{\pi }{2\sqrt{3}}\,\mathbf{I}_{1}(a_{0})[\mathfrak{%
J}_{-}(a)\mathfrak{J}_{+}^{\,\prime }(a_{0})-\mathfrak{J}_{+}(a)\mathfrak{J}%
_{-}^{\,\prime }(a_{0})]\, \\
&&-\frac{\pi \,}{2\sqrt{3}}\mathbf{I}_{2}(a_{0})[\mathfrak{J}_{-}(a)%
\mathfrak{J}_{+}(a_{0})-\mathfrak{J}_{+}(a)\mathfrak{J}_{-}(a_{0})] \\
&&+\mathfrak{J}_{+}(a)\,\Delta H_{+}(a)+\mathfrak{J}_{-}(a)\,\Delta
H_{-}(a)-Ai(-a)\,\Delta \left\{ \ln (a)\right\} ,
\end{eqnarray*}%
where%
\begin{eqnarray*}
H_{+}(a) &=&\tfrac{3^{1/3}\Gamma (2/3)^{2}}{4\pi }\,a\,\,_{1}F_{2}\tbinom{%
1/3\;;\;-a^{3}/9}{4/3,\;4/3\;}-\tfrac{\Gamma (2/3)}{\,3^{5/6}\,4}%
a^{2}\,\,_{1}F_{2}\tbinom{2/3\;;\;-a{}^{3}/9}{4/3,\;5/3\;} \\
&&+\tfrac{\sqrt{3}}{216}a_{\,\ }^{3}\,_{2}F_{3}\tbinom{1,1\;;\;-a^{3}/9}{%
2,\;2,\,\,7/3\;},
\end{eqnarray*}%
\begin{eqnarray*}
H_{-}(a) &=&-\tfrac{\pi }{3^{4/3}\Gamma (2/3)^{2}}\frac{1}{a}\,_{1}F_{2}%
\tbinom{-1/3\;;\;-a^{3}/9}{2/3,\;2/3\;}+\tfrac{\Gamma (1/3)}{3^{7/6}\,2\,}%
\,a\,\,_{1}F_{2}\tbinom{1/3\;;\;-a{}^{3}/9}{2/3,\;4/3\;} \\
&&+\tfrac{\sqrt{3}}{108\,}\,a_{\,\ }^{3}\,_{2}F_{3}\tbinom{1,1\;;\;-a^{3}/9}{%
2,\;2,\,\,5/3\;}\,,
\end{eqnarray*}%
\begin{equation*}
\Delta \,\,f\,(a)\equiv f\,(a)-f\,(a_{0}).
\end{equation*}%
We have obtained a closed form albeit complicated expression for the
Stieltjes transform $\mathbf{I}_{1}(a).$ Using the expression above with $%
a_{0}$ chosen to be the magnitude of one of the roots of $Ai^{\,\prime }(a)$
e.g. $a_{n}^{\prime },$terms containing $Ai^{\,\prime }(a_{n}^{\prime })$
vanish. Taking $a_{0}=\left| a_{1}^{\prime }\right| =1.08792997,$ and the
numerically evaluated values for $\mathbf{I}_{1}(\left| a_{1}^{\prime
}\right| )=$ $0.2109508346$ and $\mathbf{I}_{2}(\left| a_{1}^{\prime
}\right| )=$ $0.1425319307$ then allows a calculation of the integral $%
\mathfrak{I}^{\,(1)}$ as will be seen below.

\medskip The integral $\mathfrak{I}^{\,(1)}$ as represented by (3) i.e. 
\begin{equation*}
\mathfrak{I}^{\,(1)}\mathfrak{=}\frac{2}{Ai^{\prime }(0)}\sum_{n=1}^{\infty }%
\frac{1}{\left| a_{n}^{\prime }\right| }\mathbf{I}_{3}(\left| a_{n}^{\prime
}\right| ),
\end{equation*}%
when rewritten using the relation%
\begin{equation*}
\mathbf{I}_{3}(a)=\frac{1}{6}+\frac{Ai^{\,\prime }(0)}{2\,a}+\frac{Ai(0)}{%
2\,\,a^{2}}-\frac{a}{2}\mathbf{I}_{1}(a),
\end{equation*}%
gives%
\begin{equation*}
\mathfrak{\mathfrak{I}}^{\,(1)}\mathfrak{=}\frac{1}{Ai^{\,\prime }(0)}%
\sum_{n=1}^{\infty }\frac{1}{\left| a_{n}^{\prime }\right| }\{\frac{1}{3}+%
\frac{Ai^{\,\prime }(0)}{\left| a_{n}^{\prime }\right| }+\frac{Ai^{\,}(0)}{%
\left| a_{n}^{\prime }\right| ^{2}}-\left| a_{n}^{\prime }\right| \,\mathbf{I%
}_{1}(\left| a_{n}^{\prime }\right| )\}.
\end{equation*}%
\ The second and third terms in the equation above can be summed using the
Airy zeta functions $\mathcal{Z}_{k}$\ (cf. Appendix I) to give the
expression%
\begin{equation}
\mathfrak{\mathfrak{I}}^{\,(1)}\mathfrak{=}\frac{1}{Ai^{\,\prime }(0)}%
\sum_{n=1}^{\infty }\{\frac{1}{3\left| a_{n}^{\prime }\right| }-\mathbf{I}%
_{1}(\left| a_{n}^{\prime }\right| )\}.  \label{eq8}
\end{equation}%
We note that for large $a$\ (cf. Appendix II) the terms within the curly
brackets in (8) are on the order of $\left| a_{n}^{\prime }\right|
^{-2}\thicksim \left( 4/3\pi n\right) ^{4/3}$ which gives some assurance
that the sum converges.\ 

However, the value of the integral $\mathfrak{I}^{\,(1)}$ as expressed by
(8) is slowly convergent. \ Taking one hundred terms in the sum gives a
value of $-0.73273890$. \ The integral $\mathfrak{\mathfrak{I}}^{\,(1)}$
computed numerically by Maple has value $-0.81400778$. \ In contrast, the
value of $\mathfrak{\mathfrak{I}}^{\,(1)}$ using the sum in (3) which
contains one hundred terms produces $-0.81399655$, a value with error in the 
$6$ $th$ decimal place. \ We see that neither of the sums given by (3) or
(8) yield useful analytic representations for $\mathfrak{I}^{\,(1)}$. \ The
slow convergence of these series is no doubt due to contributions from terms
containing the roots $\left| a_{n}^{\prime }\right| $ with large values of $%
n.$

An approach aimed at accelerating the rate of convergence of the
representation of $\,\mathfrak{I}^{\,(1)}$ as given by equation (3) begins
by examining in f\thinspace iner detail the behavior of $\mathbf{I}_{3}(a)$
for large $a.$ \ We write%
\begin{equation*}
\mathbf{I}_{3}(a)=\int_{0}^{a}\frac{Ai(x)}{(x+a)^{3}}\,dx+\int_{a}^{\infty }%
\frac{Ai(x)}{(x+a)^{3}}\,dx.
\end{equation*}%
In the f\thinspace irst integral where $x\leq a$ the denominator can be
written as a power series in $x$ and integrated term by term. \ The
resulting value for the f\thinspace irst integral is 
\begin{eqnarray*}
\frac{1}{a}\mathbf{I}_{3},_{x\,\leq \,a}(a) &=&\frac{1}{a}\int_{0}^{a}\frac{%
Ai(x)}{(x+a)^{3}}\,dx \\
&\simeq &\frac{1}{6a^{4}}\sum_{k=0}^{n}\frac{(-1)^{k}(k+2)!}{\Gamma (k/3+1)}(%
\frac{1}{3^{1/3}a})^{k}\ldots
\end{eqnarray*}%
where we note the implied truncated series at $k=n$. \ 

In the case of the second integral i.e. where $x\geq a$ and $a$ is large, we
have%
\begin{equation*}
\frac{1}{a}\mathbf{I}_{3},_{x\geq \,\,a}\,(a)=\frac{1}{a}\int_{a}^{\infty }%
\frac{Ai(x)}{(x+a)^{3}}\,dx\simeq \frac{2}{3\sqrt{\pi }a^{9/4}}%
\int_{1}^{\infty }\frac{\exp (-\frac{2}{3}a^{3/2}x^{2})}{(1+x^{4/3})^{3}}%
\,dx,
\end{equation*}%
where the asymptotic form of $Ai(x)$ has been used. \ The largest
contribution to the latter integral occurs in the interval $1\leq x\leq 2$
so that

\begin{equation*}
\int_{a}^{\infty }\frac{Ai(x)}{(x+a)^{3}}\,dx\simeq \frac{2}{3\sqrt{\pi }%
a^{13/4}}\int_{1}^{2}\frac{\exp (-\frac{2}{3}a^{3/2}x^{2})}{(1+x^{4/3})^{3}}%
\,dx.
\end{equation*}%
Expanding the denominator in the latter integral in powers of $(x-1)$ and
integrating\ the resulting series we get%
\begin{eqnarray*}
\int_{a}^{\infty }\frac{Ai(x)}{(x+a)^{3}}\,dx &\simeq &-\tfrac{\exp (-\frac{2%
}{3}a^{3/2})}{a^{3/2}}\left( \tfrac{1345}{1728}+\tfrac{2659}{2304\,a^{3/2}}+%
\tfrac{13}{48\,a^{3}}\right) \\
&&+\tfrac{\exp (-\frac{8}{3}a^{3/2})}{a^{3/2}}\left( \tfrac{47}{108}+\tfrac{%
851}{1152\,a^{3/2}}\right) \\
&&+\tfrac{\sqrt{6\pi }}{a^{3/4}}[Erf(\tfrac{2\sqrt{6}a^{3/4}}{3})-Erf(\tfrac{%
\sqrt{6}a^{3/4}}{3})] \\
&&\cdot \left( \tfrac{1507}{5184}+\tfrac{115}{216\,a^{3/2}}+\tfrac{245}{%
1024\,a^{3}}\right) .
\end{eqnarray*}%
Using the asymptotic expression for the error function $Erf(z)$ we have upon
combining the terms in $\mathbf{I}_{3,_{\,\,x\,\,\geq \,\,a}}(a)$%
\begin{equation*}
\frac{1}{a}\mathbf{I}_{3},\,_{x\,\,\geq \,\,a}(a)\thicksim \frac{\exp (-%
\frac{2}{3}a^{\,3/2})}{\sqrt{\pi }\,a^{19/4}}.
\end{equation*}%
The exponential terms in the expression for $\mathbf{I}_{3},\,_{x\,\geq
\,\,a}(a)$ are small compared to $\mathbf{I}_{3},\,_{x\,\leq \,\,a}(a)$ and
the integration over the interval $x\geq a$ will be neglected in computing
the sums in the equation below.

With this we can write the expression for $\mathfrak{I}^{\,(1)}\mathfrak{\ }$%
using (3) approximately as%
\begin{equation*}
\mathfrak{I}^{(1)}\simeq \frac{2}{Ai^{\prime }(0)}\sum_{n=1}^{N}\frac{%
\mathbf{I}_{3}(\left| a_{n}^{\prime }\right| )}{\left| a_{n}^{\prime
}\right| }+\frac{2}{Ai^{\prime }(0)}\sum_{n=N+1}^{\infty }\frac{\mathbf{I}%
_{3},_{x\,\leq \,a}(\left| a_{n}^{\prime }\right| )}{\left| a_{n}^{\prime
}\right| },
\end{equation*}%
or introducing the inf\thinspace inite sums which contains all of the roots $%
a_{n}^{\prime }$ we have 
\begin{equation*}
\mathfrak{I}^{(1)}\simeq \frac{2}{Ai^{\prime }(0)}\sum_{n=1}^{N}\frac{1}{%
\left| a_{n}^{\prime }\right| }\left[ \mathbf{I}_{3}(\left| a_{n}^{\prime
}\right| )-\mathbf{I}_{3},_{x\,\leq \,a}(\left| a_{n}^{\prime }\right| )%
\right] +\frac{2}{Ai^{\prime }(0)}\sum_{n=1}^{\infty }\frac{\mathbf{I}%
_{3},_{x\,\leq \,a}(\left| a_{n}^{\prime }\right| )}{\left| a_{n}^{\prime
}\right| }.
\end{equation*}

In the sums in this expression, the power of $a$ in the expansion
representing $\mathbf{I}_{3},_{x\,\leq \,a}(\left| a_{n}^{\prime }\right| )$
given above whose upper limit was denoted $n,$ we note that this limit has
not been restricted and is a matter of choice as is the case for $N$. \
Summing those expressions in terms of the Airy zeta $\mathcal{Z}_{k}$ and
the incomplete Airy zeta function $\mathcal{Z}_{k}(N)$ (cf. Appendix II) we
get%
\begin{equation*}
\mathfrak{I}^{(1)}\simeq \frac{2}{Ai^{\prime }(0)}\sum_{n=1}^{N}\tfrac{1}{%
\left| a_{n}^{\prime }\right| }\mathbf{I}_{3}(\left| a_{n}^{\prime }\right|
)+\frac{1}{3Ai^{\,\prime }(0)}\sum_{k=0}^{n}\tfrac{(-1)^{k}(k+2)!}{%
3^{k/3}\Gamma (k/3+1)}\{\mathcal{Z}_{k+4}-\mathcal{Z}_{k+4}(N)\}.
\end{equation*}%
In the case where $N=10$ and $n=3,$ $\ \mathfrak{I}^{(1)}=-0.8140073597$ a
value which is accurate to seven decimal places. \ No ef\thinspace
f\thinspace ort has been made to vary $N$ and/or $n$ in an attempt to
increase the accuracy of the expression for $\mathfrak{I}^{\,(1)}.$

\subsection{Values of $\mathbf{I}_{k}(a)$ for small $a$}

Accurate values of $\mathbf{I}_{k}(a)$ for small $a$ are useful in providing
analytic expressions for the integration constants\textbf{\ }$\mathbf{I}%
_{1}(a_{0})$ and\ $\mathbf{I}_{2}(a_{0})$ i.e. the quantities which are
needed to make the solution of the dif\thinspace ferential equation (7) for $%
\mathbf{I}_{1}(a)$ unique and analytic. \ In an ef\thinspace fort to do that
we write the integral $\mathbf{I}_{n}(a)$ as

\begin{equation*}
\mathbf{I}_{\,n}(a)=\int_{a}^{\infty }\frac{\,Ai(z-a)}{z^{\,n}}\,\,dz,
\end{equation*}%
and expand the Airy function $Ai(z-a)$ in a power series in $a$,\ with the
hope that the resulting integrated series would be capable of yielding
accurate values of $\mathbf{I}_{\,n}(a)$ for $a\leq 1.$ \ The power series
for $Ai(z-a)$ i.e. 
\begin{equation*}
Ai(z-a)=\sum_{k=0}^{\infty }\frac{(-a)^{k}}{k!}\frac{d^{\,\,k}Ai(z)}{%
d\,z^{\,k}},
\end{equation*}%
is seen to require analytic expressions for the higher derivatives of the
Airy function. \ These have been studied \cite{laurpq} and are given by%
\begin{equation}
\frac{d\,^{k}Ai(z)}{d\,z^{\,k}}=\mathcal{P}_{\,k}(z)\,Ai(z)+Q\,_{k}(z)%
\,Ai^{^{\,\prime }}(z),  \label{eq9}
\end{equation}%
where $\mathcal{P}_{\,k}(z)$ and $Q\,_{k}(z)$ are polynomials. \ Recursion
relations for these polynomial have the forms 
\begin{subequations}
\begin{eqnarray}
\mathcal{P\,}_{k+2}(z) &=&z\,\mathcal{P\,}_{k}(z)+k\,\mathcal{P\,}%
_{k-1}(z)\,,  \label{eq10} \\
Q\,_{k+2}\left( z\right) &=&z\,Q\,_{k}\left( z\right) +k\,Q\,_{k-1}(z),
\label{b} \\
\mathcal{P\,}_{k+1}(z) &=&\frac{d\,\mathcal{P\,}_{k}(z)}{d\,z}+z\,Q\,_{k}(z),
\label{c} \\
Q\,_{k+1}(z) &=&\frac{d\,Q\,_{k}(z)}{d\,z}+\mathcal{P\,}_{k}(z).  \label{d}
\end{eqnarray}%
\ Initial values of $\mathcal{P}_{\,k}(z)\,$and $Q\,_{k}\left( z\right) $
i.e. $\mathcal{P}_{\,0}(z)\,=1,$ $Q\,_{0}(z)=0$ are suf\thinspace
f\thinspace icient to generate the higher polynomials.\ The f\thinspace irst
few of these are given in Table 2

\begin{center}
\begin{table}[H] \centering%
\caption{The polynomials P and Q\label{two}}%
\end{table}%

\medskip 
\begin{tabular}{|c|c|c|}
\hline
$k$ & $\mathcal{P\,}_{k}(z)$ & $Q\,_{k}\left( z\right) $ \\ \hline
$1$ & $0$ & $1$ \\ \hline
$2$ & $z$ & $0$ \\ \hline
$3$ & $1$ & $z$ \\ \hline
$4$ & $z^{2}$ & $2$ \\ \hline
$5$ & $4\,z$ & $z^{2}$ \\ \hline
$6$ & $4+z^{3}$ & $6\,z$ \\ \hline
$7$ & $9\,z^{2}$ & $10+z^{3}$ \\ \hline
$8$ & $28\,z+z^{4}$ & $12\,z^{2}$ \\ \hline
$9$ & $28+16\,z^{3}$ & $52\,z+z^{4}$ \\ \hline
$10$ & $100\,z^{2}+z^{5}$ & $80+20\,z^{3}$ \\ \hline
\end{tabular}
\end{center}

The generating functions $\xi (t,z)$ and $\lambda (t,z)$ for the polynomials 
$\mathcal{P}_{\,k}(z)$ and $Q\,_{k}\left( z\right) $ are 
\end{subequations}
\begin{eqnarray*}
\xi (t,z) &=&\sum_{k=0}^{\infty }\frac{t^{\,k}}{k!}\,\mathcal{P\,}%
_{k}(z)=\pi \lbrack Bi^{\,\prime }(z)\,Ai(z+t)-Ai^{\,\prime }(z)\,Bi(z+t)],
\\
\lambda (t,z) &=&\sum_{k=0}^{\infty }\frac{t^{\,k}}{k!}\,Q\,_{k}\left(
z\right) =\pi \lbrack Ai(z)\,Bi(z+t)-Bi(z)\,Ai(z+t)].
\end{eqnarray*}%
The expression for $Ai(z-a)$ is then given by 
\begin{equation*}
Ai(z-a)=\sum_{k=0}^{\infty }\frac{(-a)^{k}}{k!}\left\{ \mathcal{P\,}%
_{k}(z)\,Ai(z)+Q\,_{k}\left( z\right) \,Ai^{\,\prime }(z)\right\} .
\end{equation*}%
As a result, we get for the integral $\mathbf{I}_{n}(a)$ 
\begin{eqnarray}
\mathbf{I}_{n}(a) &=&\int_{a}^{\infty }\left\{ \sum_{k=0}^{\infty }\frac{%
(-a)^{k}}{k!}\mathcal{P\,}_{k}(z)\right\} \frac{\,Ai(z)}{z^{n}}\,dz
\label{eq11} \\
&&+\int_{a}^{\infty }\left\{ \sum_{k=0}^{\infty }\frac{(-a)^{k}}{k!}%
Q\,_{k}\left( z\right) \,\right\} \frac{Ai^{\,\prime }(z)}{z^{n}}\,dz. 
\notag
\end{eqnarray}%
or in terms of the generating functions%
\begin{equation}
\mathbf{I}_{n}(a)=\int_{a}^{\infty }\left\{ \xi (-a,z)\,Ai(z)\,+\lambda
(-a,z)\,Ai^{\,\prime }(z)\right\} \frac{\,dz}{z^{n}}.  \label{eq12}
\end{equation}

Initially we choose to use the former expression (11) as a means of
computing $\mathbf{I}_{3}(a)$ i.e. 
\begin{equation}
\mathbf{I}_{3}(a)=\sum_{k=0}^{\infty }\frac{(-a)^{k}}{k!}\left\{
\int_{a}^{\infty }\frac{\mathcal{P\,}_{k}(z)\,Ai(z)}{z^{3}}%
\,dz+\int_{a}^{\infty }\frac{Q\,_{k}\left( z\right) \,Ai^{\,\prime }(z)}{%
z^{3}}\,dz\right\} .  \label{eq13}
\end{equation}%
The integral containing $Ai^{\,\prime }(z)$ in (13) can be integrated by
parts with the result%
\begin{eqnarray*}
\int_{a}^{\infty }\frac{Q\,_{k}\left( z\right) \,Ai^{\,\prime }(z)}{z^{3}}%
\,dz &=&-\frac{Q\,_{k}\left( a\right) \,Ai(a)}{a^{3}} \\
&&+\int_{a}^{\infty }\left[ \frac{3\,Q\,_{k}\left( z\right) \,}{z^{4}}+\frac{%
\mathcal{P\,}_{k}(z)}{z^{3}}-\frac{Q\,_{k+1}\left( z\right) }{z^{3}}\right]
\,Ai(z)\,dz
\end{eqnarray*}%
Combining the integrals in (13) we get%
\begin{equation}
\mathbf{I}_{3}(a)=\sum_{k=0}^{\infty }\frac{(-a)^{k}}{k!}\int_{a}^{\infty
}F_{k}(z)\,Ai(z)\,dz-\frac{Ai(a)}{a^{3}}\sum_{k=0}^{\infty }\frac{(-a)^{k}}{%
k!}Q\,_{k}\left( a\right) ,  \label{eq14}
\end{equation}%
where 
\begin{equation*}
z^{4}F_{k}(z)=z\{2\,\mathcal{P\,}_{k}(z)\,-Q\,_{k+1}\left( z\right)
\}+3\,Q\,_{k}\left( z\right) .
\end{equation*}%
\ The sum in (14) which contains the polynomials $Q\,_{k}\left( a\right) $
can be related to the Airy functions using its generating function. We get
the closed form expression 
\begin{equation*}
\sum_{k=0}^{\infty }\frac{(-a)^{k}}{k!}\,Q\,_{k}\left( a\right) =\pi Ai(0)[%
\sqrt{3}Ai(a)-Bi(a)]=\pi Ai(0)\mathfrak{I}_{-}(-a).
\end{equation*}%
\ The degrees of the polynomials in $z^{4}F_{k}(z)$ as seen below in the
matrix $\mathfrak{A}$\ are a complicated function of the index $k$. \ As a
result, rearranging the integrations in (14) in terms of increasing powers
of $z$ is in general also complicated. \ We write the integrand\thinspace $%
F_{k}(z)$ as%
\begin{equation*}
F_{k}(z)=\frac{1}{z^{4}}\sum_{i=0}^{N_{k}}\mathfrak{A}_{k,\,i}(a)\,\,z^{\,i},
\end{equation*}%
and note that in the case of even and odd values of $k$ the upper limits $%
N_{k}$ in the sums are just 
\begin{eqnarray*}
F_{2k}(z) &=&\frac{1}{z^{4}}\sum_{i=0}^{k+1}\mathfrak{A}_{2k,\,i}(a)\,\,z^{%
\,i}, \\
F_{2k+1}(z) &=&\frac{1}{z^{4}}\sum_{i=0}^{k}\mathfrak{A}_{2k+1,\,i}(a)\,%
\,z^{\,i},
\end{eqnarray*}%
as can be seen by an inspection of the elements $\mathfrak{A}_{k}\mathfrak{,}%
_{i}(a)$ which are displayed in the matrix $\mathfrak{A}$

\begin{equation*}
\mathfrak{A}=\left[ 
\begin{array}{cccccccccc}
0 & 1 & 0 & 0 & 0 & 0 & 0 & 0 & 0 & \cdots \\ 
3 & 0 & 0 & 0 & 0 & 0 & 0 & 0 & 0 & \cdots \\ 
0 & 0 & 1 & 0 & 0 & 0 & 0 & 0 & 0 & \cdots \\ 
0 & 3 & 0 & 0 & 0 & 0 & 0 & 0 & 0 & \cdots \\ 
6 & 0 & 0 & 1 & 0 & 0 & 0 & 0 & 0 & \cdots \\ 
0 & 0 & 5 & 0 & 0 & 0 & 0 & 0 & 0 & \cdots \\ 
0 & 16 & 0 & 0 & 1 & 0 & 0 & 0 & 0 & \cdots \\ 
30 & 0 & 0 & 9 & 0 & 0 & 0 & 0 & 0 & \cdots \\ 
0 & 0 & 40 & 0 & 0 & 1 & 0 & 0 & 0 & \cdots \\ 
0 & 132 & 0 & 0 & 15 & 0 & 0 & 0 & 0 & \cdots \\ 
240 & 0 & 0 & 100 & 0 & 0 & 1 & 0 & 0 & \cdots \\ 
0 & 0 & 440 & 0 & 0 & 23 & 0 & 0 & 0 & \cdots \\ 
0 & 1480 & 0 & 0 & 230 & 0 & 0 & 1 & 0 & \cdots \\ 
\cdots & \cdots & \cdots & \cdots & \cdots & \cdots & \cdots & \cdots & 
\cdots & \cdots%
\end{array}%
\right]
\end{equation*}%
The expression for the integral $\mathbf{I}_{3}(a)$ is then 
\begin{eqnarray*}
\mathbf{I}_{3}(a) &=&\pi Ai(0)\mathfrak{I}_{-}(-a)-3\pi Ai(0)\mathfrak{I}%
_{-}(a)\int_{a}^{\infty }\tfrac{Ai(z)}{z^{4}}dz \\
&&+\sum_{i=1}^{\infty }a^{2i}\{\sum_{k=0}^{\infty }\tfrac{(-a^{k})}{(2k+2i)!}%
[\mathfrak{A}_{2k+2i,i}-\tfrac{a}{(2k+2i+1)}\mathfrak{A}_{2k+2i+1,i}]\}%
\int_{a}^{\infty }\tfrac{Ai(z)}{z^{\,4-i}}dz
\end{eqnarray*}%
where we have used the sum 
\begin{eqnarray*}
\sum_{k=0}^{\infty }\frac{(-a)^{k}}{k!}\mathfrak{A}_{k,0} &=&\frac{\sqrt{3}}{%
2Ai^{\,\prime }(0)}[\sqrt{3}Ai(-a)-Bi(-a)] \\
&=&-3\pi Ai(0)\mathfrak{I}_{-}(a).
\end{eqnarray*}%
The latter sum follows from an inspection of the elements of $\mathfrak{A}$
i.e. those occurring in its f\thinspace irst column i.e. 
\begin{eqnarray*}
\mathfrak{A}_{3k+1,0} &=&3^{k+1}\Gamma (k+2/3)/\Gamma (2/3), \\
\mathfrak{A}_{3k,0} &=&0, \\
\mathfrak{A}_{3k+2,0} &=&0.
\end{eqnarray*}%
Calculation of the integrals occurring in the expression for $\mathbf{I}%
_{3}(a)$ appearing above in (14) require analytical expressions for the
integrals $\int_{a}^{\infty }z^{n}Ai(z)\,dz$ and $\int_{a}^{\infty
}z^{n}Ai^{\,\prime }(z)\,dz.$ \ These are given below.

\subsection{The incomplete Mellin transforms of $\ Ai(z)$ and $Ai^{\,\prime
}(z)\,$}

\bigskip

We def\thinspace ine the integrals ( $a>0$) 
\begin{eqnarray*}
I_{n}(a) &=&\int_{a}^{\infty }\,z^{n}Ai\,(z)\,dz, \\
I_{n}^{\,\,\prime }(a) &=&\int_{a}^{\infty }\,z^{n}Ai\,^{\prime }(z)\,dz,
\end{eqnarray*}%
and note that integration by parts and use of (2) yields the third-order
recurrence relation for all $n$%
\begin{equation}
I_{n}(a)=(n-1)(n-2)\,I_{n\,-3}(a)-a^{\,n-1}Ai^{\,\prime
}(a)+(n-1)\,a^{n-2}Ai(a).  \label{eq15}
\end{equation}%
Initial values of $I_{n}(a)$ i.e. $I_{0}(a),\,I_{-1}(a),$ and $I_{-2}(a)$
are irreducible and are required to obtain the general solution to this
dif\thinspace f\thinspace erence equation. \ For $I_{0}(a)$ we have
immediately 
\begin{equation*}
I_{0}(a)=\pi \left[ Ai(a)\,Gi^{\,\prime }(a)-Ai^{\,\prime }(a)\,Gi(a)\right]
,
\end{equation*}%
where $Gi$ and $Gi^{\,\prime }$ are the Scorer functions \cite{nbsscorer}.
The integral $I_{0}(a)$ has also been denoted by $Ai_{1}(a)$ and was studied
by Aspnes \cite{Aspnes}. An expression for $I_{0}(a)$ in terms of the more
accessible hypergeometric functions being 
\begin{equation*}
I_{0}(a)=\tfrac{1}{3}-\tfrac{a}{3^{2/3}\Gamma (2/3)}\,_{1}F_{2}\tbinom{%
1/3\;;\;a^{3}/9}{2/3,\;4/3\;}+\tfrac{3^{1/6}a^{2}\Gamma (2/3)}{4\pi }%
\,_{1}F_{2}\tbinom{2/3\;;\;a^{3}/9}{4/3,\;5/3\;}.\,
\end{equation*}%
The case for $I_{-1}(a)$ is more complicated and Maple gives ($\gamma $ is
Euler's constant) 
\begin{eqnarray*}
I_{-1}(a) &=&\frac{3^{1/6}a\,\Gamma (2/3)}{2\pi }\,_{1}F_{2}\tbinom{%
1/3\;;\;a^{3}/9}{4/3,\;4/3\;}-\frac{3^{1/3}a^{3}\,}{54\,\Gamma (2/3)}%
\,_{2}F_{3}\tbinom{1,\;1\,;\;a^{3}/9}{2,\;2,\,5/3} \\
&&-\frac{3^{1/3}}{18\,\Gamma (2/3)}\left\{ \ln (a^{6}/3)+4\gamma +\sqrt{3}%
\pi /3\right\} .
\end{eqnarray*}%
Similarly, the value of the integral $I_{-2}(a)$ is given by Maple as%
\begin{eqnarray*}
I_{-2}(a) &=&\frac{1}{3^{2/3}\,a\,\Gamma (2/3)}\,_{1}F_{2}\tbinom{%
-1/3\;;\;a^{3}/9}{2/3,\;2/3\;}+\frac{3^{1/6}\,\Gamma \,(2/3)\,a^{3}\,}{72\pi 
}\,_{2}F_{3}\tbinom{1,1\;;\;a^{3}/9}{2,2,7/3\;} \\
&&+\frac{3^{1/6}\Gamma (2/3)}{12\pi }\left\{ \ln (a^{6}/3)+4\gamma -6-\sqrt{3%
}\pi /3\right\} \text{ }.
\end{eqnarray*}%
\ The values of the f\thinspace irst few of the required $I_{n}(a)$
integrals obtained using the recurrence relations (15) are%
\begin{eqnarray*}
I_{-6}(a) &=&\tfrac{1}{40}[I_{0}(a)+(1/a+2/a^{4})\,Ai^{\,\prime
}(a)+(1/a^{2}+8/a^{5})\,Ai(a)], \\
I_{-5}(a) &=&\tfrac{1}{12}[I_{-2}(a)+Ai^{\,\prime }(a)/a^{3}+3Ai(a)/a^{4}],
\\
I_{-4}(a) &=&\tfrac{1}{6}\left[ I_{-1}(a)+Ai^{\,\prime
}(a)/a^{2}+2Ai(a)/a^{3}\right] , \\
I_{-3}(a) &=&\tfrac{1}{2}\left[ I_{0}(a)+Ai^{\,\prime }(a)/a+Ai(a)/a^{2}%
\right] , \\
I_{1}(a) &=&-Ai^{\,\prime }(a), \\
I_{2}(a) &=&-a\,Ai^{\,\prime }(a)+Ai(a), \\
I_{3}(a) &=&2\,I_{0}(a)-a^{2}Ai^{\,\prime }(a)+2a\,Ai(a), \\
I_{4}(a) &=&-(6+a^{3})\,Ai^{\,\prime }(a)+3a^{2}Ai(a), \\
I_{5}(a) &=&-(12\,a+a^{4})\,Ai^{\,\prime }(a)+(12+4\,a^{3})\,Ai(a), \\
I_{6}(a) &=&40\,I_{0}(a)-(20\,a^{2}+a^{5})\,Ai^{\,\prime
}(a)+(40\,a+5\,a^{4})\,Ai(a).
\end{eqnarray*}%
In general the $I_{n}(a)$ integrals for $n\geq 0$ are given by (cf. Appendix
III) 
\begin{eqnarray*}
I_{3k}(a) &=&\frac{(3k)!}{3^{k}\,k\,!}\left\{ 
\begin{array}{c}
I_{0}(a)+a\,Ai(a)\sum_{l=0}^{k-1}\tfrac{(3a^{3})^{l}\,l!}{(3l+1)!}- \\ 
a^{2}Ai^{\,\prime }(a)\sum_{l=0}^{k-1}\tfrac{(3a^{3})^{l}\,l!}{(3l+2)!}%
\end{array}%
\right\} , \\
I_{3k+1}(a) &=&3^{2k}k!\Gamma (k+2/3)\left\{ 
\begin{array}{c}
\tfrac{a^{2}}{3}Ai^{\,}(a)\sum_{l=0}^{k-1}\tfrac{(a^{3}/9)^{l}\,}{%
l!\,\,\Gamma (l+5/3)} \\ 
-Ai^{\,\prime }(a)\sum_{l=0}^{k}\tfrac{(a^{3}/9)^{l}\,}{l!\,\,\Gamma (l+2/3)}%
\end{array}%
\right\} , \\
I_{3k+2}(a) &=&3^{2k}k!\Gamma (k+4/3)\left\{ 
\begin{array}{c}
3Ai^{\,}(a)\sum_{l=0}^{k}\tfrac{(a^{3}/9)^{l}\,}{l!\,\,\,\Gamma (l+1/3)} \\ 
-a\,Ai^{\,\prime }(a)\sum_{l=0}^{k}\tfrac{(a^{3}/9)^{l}\,}{l!\,\,\Gamma
(l+4/3)}%
\end{array}%
\right\} ,
\end{eqnarray*}%
and $\ $the\ $I_{-n}(a)$ for $n\geq 0$ are given by%
\begin{eqnarray*}
I_{-3k}(a) &=&\tfrac{3^{k}\,k\,!}{(3k)!}\left\{ 
\begin{array}{c}
I_{0}(a)+\frac{Ai(a)}{a^{2}}\sum_{l=0}^{k-1}\tfrac{(3l+1)!}{(3a^{3})^{l}\,l!}
\\ 
+\frac{Ai^{\,\prime }(a)}{a}\sum_{l=0}^{k-1}\tfrac{(3l+2)!}{(3a^{3})^{l}\,l!}%
\end{array}%
\right\} , \\
&& \\
I_{-(3k+1)}(a) &=&\tfrac{1}{3^{2k}k!\,\Gamma (k+2/3)}\left\{ 
\begin{array}{c}
\Gamma (2/3)\,I_{-1}(a) \\ 
+\tfrac{3}{a^{3}}Ai^{\,}(a)\sum_{l=0}^{k-1}\tfrac{l\,!\,\,\Gamma (l+5/3)\,}{%
(a^{3}/9)^{l}} \\ 
+\tfrac{1}{a^{2}}Ai^{\,\prime }(a)\sum_{l=0}^{k-1}\tfrac{l\,!\,\,\Gamma
(l+2/3)\,}{(a^{3}/9)^{l}}%
\end{array}%
\right\} , \\
&& \\
I_{-(3k+2)}(a) &=&\tfrac{1}{3^{2k}k!\,\Gamma (k+4/3)}\left\{ 
\begin{array}{c}
\Gamma (4/3)\,I_{-2}(a) \\ 
+\tfrac{3}{a^{4}}Ai^{\,}(a)\sum_{l=0}^{k-1}\tfrac{(l+1)!\,\,\Gamma (l+4/3)\,%
}{(a^{3}/9)^{l}} \\ 
+\tfrac{1}{a^{3}}Ai^{\,\prime }(a)\sum_{l=0}^{k-1}\tfrac{l!\,\,\Gamma
(l+4/3)\,}{(a^{3}/9)^{l}}%
\end{array}%
\right\} .
\end{eqnarray*}%
\ In cases where the upper limits of these sum contain negative values the
sums above are empty. \ 

\bigskip

The primed integrals $I_{n}^{^{\;\prime }}(a)$ when integrated by parts
gives for $n\geq 0$ 
\begin{equation*}
I_{n}^{^{\;\prime }}(a)=-[a^{\,n}Ai(a)+n\,I_{n-1}(a)].
\end{equation*}%
In addition we also have the useful relation $(n\geq 2)$ 
\begin{equation*}
I_{-n}^{^{\;\,\prime }}(a)=[Ai^{\,\prime
}(a)/a^{\,n-1}+\,I_{-n+2\,}(a)]/(n-1),\hspace{0.25in}
\end{equation*}

General expressions for the integrals $I_{n}^{^{\;\prime }}(a)$ are given by%
\begin{eqnarray*}
I_{3k}^{\,\,\prime }(a) &=&9^{k}k!\Gamma (k+1/3)\{-Ai(a)\sum_{l=0}^{k}\tfrac{%
(a^{3}/9)^{l}\,}{l!\,\Gamma (l+1/3)} \\
&&+\tfrac{1}{2}a\,Ai^{\,\prime }(a)\sum_{l=0}^{k-1}\tfrac{(a^{3}/9)^{l}\,}{%
l!\Gamma (l+4/3)}\}, \\
I_{3k+1}^{\,\,\prime }(a) &=&\tfrac{(3k+1)!}{3^{k}}\{-I_{0}(a)-a\,Ai(a)%
\sum_{l=0}^{k}\tfrac{(3a^{3})^{l}\,l!\,}{(3k+1)!} \\
&&+a^{2}\,Ai^{\,\prime }(a)\sum_{l=0}^{k-1}\tfrac{(3a^{3})^{l}\,\,l!}{(3k+2)!%
}\}, \\
I_{3k+2}^{\,\,\prime }(a) &=&9^{k}k!\Gamma (k+5/3)\{-a^{2}Ai(a)\sum_{l=0}^{k}%
\tfrac{(a^{3}/9)^{l}\,}{l!\,\Gamma (l+5/3)} \\
&&+3Ai^{\,\prime }(a)\sum_{l=0}^{k-1}\tfrac{(a^{3}/9)^{l}\,}{l!\Gamma (l+2/3)%
}\}.
\end{eqnarray*}%
Some of the $I_{n}^{^{\;\prime }}(a)$ integrals\ appearing in $\mathbf{I}%
_{3}(a)$ are explicitly given by

\begin{eqnarray*}
I_{-4}^{\,\,\,\prime }(a) &=&I_{-2}(a)/3-Ai^{\,\,\prime }(a)/3a^{3}, \\
I_{-3}^{\,\,\,^{\prime }}(a) &=&I\,_{-1}(a)/2-Ai^{\,\,\prime }(a)/2a^{2}, \\
I_{-2}^{\,\,\,\prime }(a) &=&I\,_{0}(a)+Ai^{\,\,\prime }(a)/a, \\
I_{-1}^{\,\,\prime }(a) &=&I_{-2}(a)-Ai(a)/a, \\
I_{0}^{\,\,\prime }(a) &=&-Ai(a), \\
I_{1}^{\,\,\prime }(a) &=&-I_{0}(a)-a\,Ai(a), \\
I_{2}^{\,\,\prime }(a) &=&-a^{2}Ai(a)+2Ai^{\,\,\prime }(a), \\
I_{3}^{\,\,\prime }(a) &=&-(3+a^{3})Ai(a)+3a\,Ai^{\,\,\prime }(a), \\
I_{4}^{\,\,\prime }(a) &=&-8I_{0}(a)-(8a+a^{4})Ai(a)+4a^{4}Ai^{\,\,\prime
}(a), \\
I_{5}^{\,\,\prime }(a) &=&-(15a^{2}+a^{5})Ai(a)+5(6+a^{3})Ai^{\,\,\prime
}(a), \\
I_{6}^{\,\,\prime }(a)
&=&-(72+24a^{3}+a^{6})Ai(a)+6a\,(12+a^{3})Ai^{\,\,\prime }(a), \\
&&\cdots
\end{eqnarray*}%
(See Appendix III for an alternate means of computing the integrals $%
I_{n}(a) $ and $I_{n}^{\,\prime }(a)$.)

In terms of these integrals we have for $\mathbf{I}_{3}(a)$%
\begin{eqnarray*}
\mathbf{I}_{3}(a) &=&\pi Ai(0)[\sqrt{3}Ai(a)-Bi(a)]+\tfrac{\sqrt{3}}{%
2Ai^{\,\prime }(0)}[\sqrt{3}Ai(-a)-Bi(-a)]\,I_{-4}(a) \\
&&+\sum_{i=1}^{\infty }a^{2i}\{\sum_{k=0}^{\infty }\tfrac{(-a^{k})}{(2k+2i)!}%
[\mathfrak{A}_{2k+2i,i}-\tfrac{a}{(2k+2i+1)}\mathfrak{A}_{2k+2i+1,i}]%
\}I_{i-4}(a).
\end{eqnarray*}%
Using this expression with inf\thinspace inite sums truncated to include
powers of $a$ \ up to $a^{10}$ we obtain a value of $\mathbf{I}_{3}(\left|
a_{1}^{\prime }\right| )$ i.e.%
\begin{equation*}
\mathbf{I}_{3}(\left| a_{1}^{\prime }\right| )=0.1045962658,
\end{equation*}%
whereas the numerical value of $\mathbf{I}_{3}(\left| a_{1}^{\prime }\right|
)=0.1045955174.$ Using the relation connecting $\mathbf{I}_{3}(a)$ and $%
\mathbf{I}_{1}(a)$ we have 
\begin{equation*}
\mathbf{I}_{1}(\left| a_{1}^{\prime }\right| )=0.2082343827,
\end{equation*}%
as compared with the numerical value of $\mathbf{I}_{1}(\left| a_{1}^{\prime
}\right| )$ i.e. $0.2082347508$ a value accurate to six decimal places. \
Since it is possible to increase the accuracy of these calculations as will
be seen below we postpone the calculation of $\mathbf{I}_{2}(\left|
a_{1}^{\prime }\right| )$. \ 

\bigskip

In obtaining the results above we have used closed form expression for the
leading terms of some of the inf\thinspace inite sums appearing in equation
(11). \ Now we consider the possibility of continuing the process of
including higher order contributions due to $a$ in a more systematic way.
Expanding the generating functions of\ equation (12) in power series in $z$
we obtain an expression which can be readily integrated i.e.%
\begin{equation*}
\mathbf{I}_{n}(a)=\sum_{i=0}^{\infty }\tfrac{\xi ^{(i)}(-a,0)}{i!}%
\int_{a}^{\infty }\,z^{i-n}Ai(z)\,dz+\sum_{i=0}^{\infty }\tfrac{\lambda
^{(i)}(-a,0)}{i!}\int_{a}^{\infty }z^{i-n}Ai^{\,\prime }(z)\,dz,
\end{equation*}%
or in terms of the incomplete Mellin transforms $I_{i}(a)$ and $%
I_{i}^{\,\,\prime }(a)$ we have%
\begin{equation*}
\mathbf{I}_{n}(a)=\sum_{i=0}^{\infty }\tfrac{\xi ^{(i)}(-a,0)}{i!}%
I_{i-n}(a)+\sum_{i=0}^{\infty }\tfrac{\lambda ^{(i)}(-a,0)}{i!}%
I_{i-n}^{\,\,\prime }(a).
\end{equation*}%
Since the integrals $I_{n}^{\,\,}(a)$ and $I_{n}^{\,\,\prime }(a)$ are
related we have%
\begin{eqnarray*}
\mathbf{I}_{n}(a) &=&\sum_{i=1}^{\infty }\frac{\xi
^{(n+i)}(-a,0)I_{i}(a)-\lambda ^{(n+i)}(-a,0)[a^{i}Ai(a)+i\,I_{i-1}^{\,}(a)]%
}{(n+i)!} \\
&&+\sum_{i=0}^{n}\frac{\xi ^{(n-i)}(-a,0)I_{-\,i}(a)-\lambda
^{(n-i)}\,(-a,0)[Ai(a)/a^{i}-iI_{-i-1}(a)]}{(n-i)!},
\end{eqnarray*}%
where the sums containing integral transforms with positive and negative
indexes have been written as separate sums. \ The f\thinspace irst few
coef\thinspace f\thinspace icients $\xi ^{(i)}(-a,0)$ and $\lambda
^{(i)}(-a,0)$ (hereafter denoted by $\xi ^{(i)}$ and $\lambda ^{(i)}$) are%
\begin{eqnarray}
\xi ^{(0)} &=&-\pi Ai^{\,\prime }(0)\mathfrak{J}_{+}(a),  \notag \\
\xi ^{\left( 1\right) } &=&-\pi Ai^{\,\prime }(0)\mathfrak{J}_{+}^{\,\prime
}(a),  \notag \\
\xi ^{\left( 2\right) } &=&\pi Ai(0)\mathfrak{J}_{-}(a)+a\,\,\pi
Ai^{\,\prime }(0)\mathfrak{J}_{+}(a),  \notag \\
\xi ^{\left( 3\right) } &=&3\pi Ai(0)\mathfrak{J}_{-}^{\,\prime }(a)-\pi
Ai^{\,\prime }(0)[3\,\mathfrak{J}_{+}(a)-a\,\mathfrak{J}_{+}^{\,\prime }(a)],
\label{eq16}
\end{eqnarray}%
\begin{eqnarray}
\lambda ^{(0)} &=&-\pi Ai(0)\mathfrak{J}_{-}(a),  \notag \\
\lambda ^{(1)} &=&-\pi Ai(0)\mathfrak{J}_{-}^{\,\prime }(a)+\pi Ai^{\,\prime
}(0)\mathfrak{J}_{+}(a),  \notag \\
\lambda ^{(2)} &=&a\,\pi Ai(0)\mathfrak{J}_{-}(a)+2\pi Ai^{\,\prime }(0)%
\mathfrak{J}_{+}^{\,\prime }(a),  \notag \\
\lambda ^{(3)} &=&-\pi Ai(0)[2\,\mathfrak{J}_{-}(a)-a\,\mathfrak{J}%
_{-}^{\,\prime }(a)]-3a\,\pi Ai^{\,\prime }(0)\mathfrak{J}_{+}(a).
\label{eq17}
\end{eqnarray}%
We get for $\mathbf{I}_{n}(a)$%
\begin{equation*}
\mathbf{I}_{n}(a)=-\frac{Ai(a)}{a^{n}}\sum_{i=0}^{\infty }\frac{%
a^{\,i}\lambda ^{(i)}}{i!}+\sum_{i=0}^{\infty
}b_{i}\,I_{i}(a)+\sum_{i=1}^{n+1}b_{-i\,}I_{-\,i}(a).
\end{equation*}%
where the coef\thinspace f\thinspace icients $b_{i}$ are given by%
\begin{equation*}
b_{i}=\frac{1}{(n+1+i)!}\left\{ (n+1+i)\xi ^{(n+i)}-(1+i)\lambda
^{(n+1+i)}\right\} .
\end{equation*}%
A general expression for $\mathbf{I}_{n}(a)$ for all $n$ is quite
complicated and is given in Appendix IV below. \ The corresponding
expression for $\mathbf{I}_{3}(a)$ being 
\begin{eqnarray*}
\mathbf{I}_{3}(a) &=&I_{0}(a)\{\frac{2\pi }{\sqrt{3}}\omega
_{1,1}+\sum_{k=0}^{\infty }\Omega _{k,0}(3,a)\} \\
&&+I_{-1}(a)\Gamma (2/3)\left\{ \omega _{0,2}+\omega _{1,2}\right\}
+I_{-2}(a)\frac{1}{3}\Gamma (1/3)\omega _{0,3} \\
&&+Ai(a)\left\{ 
\begin{array}{c}
-\frac{1}{a^{3}}\pi Ai(0)\mathfrak{J}_{-}(-a)+\omega _{1,1}\widehat{\sigma }%
_{1,1}+\omega _{1,2}\widehat{\sigma }_{1,2} \\ 
+\sum_{\mu =0}^{2}\sum_{k=0}^{\infty }\Omega _{k,\mu }(3,a)\sigma _{k,\mu }%
\end{array}%
\right\} \\
&&+Ai^{\,\prime }(a)\left\{ \omega _{1,1}\widehat{\sigma }_{1,1}^{\,\prime
}+\omega _{1,2}\widehat{\sigma }_{1,2}^{\,\prime }-\sum_{\mu
=0}^{2}\sum_{k=0}^{\infty }\Omega _{k,\mu }(3,a)\sigma _{k,\mu }^{\,\prime
}\right\} ,
\end{eqnarray*}%
where the various terms appearing above are def\thinspace ined in Appendix
IV. \ Taking ten terms in the sums over $k$ together with $a=\left|
a_{1}^{\prime }\right| =1.0187929716$ the terms involving $Ai^{\,\prime
}(\left| a_{1}^{\prime }\right| )$ vanish and we obtain $\mathbf{I}%
_{3}(\left| a_{1}^{\prime }\right| )=0.1045955174$ a value accurate to 
\textit{ten decimal} places. In a similar calculation for $\mathbf{I}%
_{4}(1.0187929716)$ a value of $0.08085800094$ was obtained which is also
accurate to \textit{ten decimals}. \ From these expressions the required
expressions for $\mathbf{I}_{1}(\left| a_{1}^{\prime }\right| ),$ and $%
\mathbf{I}_{2}(\left| a_{1}^{\prime }\right| )$ can be obtained using the
relations%
\begin{eqnarray*}
\mathbf{I}_{1}(a) &=&\frac{1}{3a}+\frac{Ai^{\,\prime }(0)}{a^{2}}+\frac{%
Ai^{\,}(0)}{a^{3}}-\frac{2}{a}\mathbf{I}_{3}(a), \\
\mathbf{I}_{2}(a) &=&\frac{1}{3a^{2}}+\frac{2Ai^{\,\prime }(0)}{a^{3}}+\frac{%
3Ai^{\,}(0)}{a^{4}}-\frac{2}{a^{3}}\mathbf{I}_{3}(a)-\frac{6}{a}\mathbf{I}%
_{4}(a).
\end{eqnarray*}

\section{An integral containing $Ai^{\,\prime }(x)^{2}\ln (Ai^{\,\prime
}(x)) $}

We next consider the integral 
\begin{equation*}
\mathfrak{I}^{\left( 2\right) }\mathfrak{=}\int_{0}^{\infty }\left( \frac{%
Ai^{\,\prime }(x)}{Ai^{\,\prime }(0)}\right) ^{2}\ln \left( \frac{%
Ai^{\,\prime }(x)}{Ai^{\,\prime }(0)}\right) \,dx,
\end{equation*}%
with a value of $-\,0.2636317105$ when computed numerically. \ We have%
\begin{equation*}
\mathfrak{I}^{(2)}\mathfrak{=}\frac{1}{3Ai^{\,\prime }(0)^{2}}%
\sum_{k=1}^{\infty }\mathbf{J}(\left| a_{k}^{\,\prime }\right| ),
\end{equation*}%
with%
\begin{equation*}
\mathbf{J}(a)=\frac{1}{a}\int_{0}^{\infty }\frac{x}{\left( x+a\right) }%
\{2\,Ai(x)\,Ai^{\,\prime }(x)+x\,Ai^{\,\prime }(x)^{2}-x^{2}\,Ai(x)^{2}\}dx,
\end{equation*}%
where we have once again used the Weierstrass inf\thinspace inite product
representation for the derivative of the $\ln \left( Ai^{\,\prime
}(x)/Ai^{\,\prime }(0)\right) $ and have integrated by parts.\ The
inf\thinspace inite series expression for $\mathfrak{I}^{\left( 2\right) }$
like $\mathfrak{I}^{\,\left( 1\right) }$ is slowly convergent. \ Taking 50
terms in the series above yields a value of $-0.2343590038$ a poor estimate
of the integral in question. \ As in the case of $\mathfrak{I}^{\,\left(
1\right) },$ the rate of convergence of the sum representing $\mathfrak{I}%
^{\,\left( 2\right) }$ can be accelerated as will be seen below.

\subsection{Generalized Stieltjes transforms of $Ai(x)^{2}$ and $%
Ai^{\,\prime }(x)^{2}$}

In order to proceed further with the analysis of $\mathfrak{I}^{\,\left(
2\right) },$ it is useful to def\thinspace ine the integrals%
\begin{equation*}
J_{n}(a)=\int_{0}^{\infty }\frac{Ai(x)^{2}}{\left( x+a\right) ^{n}}\,dx,%
\text{ \ and \ }J_{n}^{\,\,\prime }(a)=\int_{0}^{\infty }\frac{Ai^{\,\prime
}(x)^{2}}{\left( x+a\right) ^{n}}\,dx,
\end{equation*}%
and to examine their properties. \ Integration of $J_{n}(a)$ by parts gives%
\begin{equation}
(n-1)\,J_{n}(a)=-\frac{Ai^{\,\prime }(0)^{2}}{a^{\,n}}+a\,n\,J_{n+1}(a)+n%
\,J_{n+1}^{\,\,\prime }(a),  \label{eq18}
\end{equation}%
and integrating 
\begin{equation*}
\int_{0}^{\infty }\frac{x\,Ai(x)^{2}}{(x+a)^{n}}\,dx,
\end{equation*}%
by parts gives 
\begin{equation}
J_{n}^{\,\prime }(a)=-\tfrac{1}{a^{n}}Ai(0)\,Ai^{\,\prime }(0)-\tfrac{n}{%
2a^{n+1}}Ai(0)^{2}-J_{n-1}(a)+a\,J_{n}(a)+\tfrac{n(n+1)}{2}J_{n+2}(a).
\label{eq19}
\end{equation}%
In terms of these quantities the expression for $\mathfrak{I\,}^{\left(
2\right) }$ becomes 
\begin{equation*}
\mathfrak{I}^{\,\left( 2\right) }\,\mathfrak{=}\frac{1}{3\,Ai^{\,\prime
}(0)^{2}}\sum_{k=1}^{\infty }\left\{ 
\begin{array}{c}
\tfrac{1}{10\left| a_{k}^{\prime }\right| }Ai(0)^{2}+\tfrac{1}{3}%
Ai(0)Ai^{\,\prime }(0)-\left| a_{k}^{\prime }\right| \,Ai^{\,\prime }(0)^{2}
\\ 
+\left| a_{k}^{\prime }\right| \,J_{1}^{\,\,\prime }(\left| a_{k}^{\prime
}\right| )+\left| a_{k}^{\prime }\right| ^{2}J_{1}(\left| a_{k}^{\prime
}\right| )-J_{2}(\left| a_{k}^{\prime }\right| )%
\end{array}%
\right\} ,
\end{equation*}%
Using (18, 19) we obtain a third-order dif\thinspace ference equation for
the integrals $J_{n}(a)$ i.e. 
\begin{eqnarray*}
&&(2n-1)\,J_{n}(a)-2a\,n\,J_{n+1}(a)-\tfrac{1}{2}n(n+1)(n+2)\,J_{n+3}(a) \\
&=&-\tfrac{1}{a^{n}}Ai^{\,\prime }(0)^{2}-\tfrac{n}{a^{n+1}}%
Ai(0)\,Ai^{\,\prime }(0)-\tfrac{n(n+1)}{2\,a^{n+2}}Ai(0)^{2},
\end{eqnarray*}%
which can be rewritten as the dif\thinspace ferential equation i.e.%
\begin{eqnarray}
&&\tfrac{1}{2}\frac{d^{\,\,3}J_{n}(a)}{d\,a^{\,3}}+2a\frac{d\,J_{n}(a)}{d\,a}%
+(2n-1)\,J_{n}(a)=  \notag \\
&&-\tfrac{1}{a^{n}}Ai^{\,\prime }(0)^{2}-\tfrac{n}{a^{n+1}}%
Ai(0)\,Ai^{\,\prime }(0)-\tfrac{n\,(n+1)}{2\,a^{n+2}}Ai(0)^{2},
\label{eq20n}
\end{eqnarray}%
having used the f\thinspace irst and third of the dif\thinspace ferential
relations%
\begin{equation*}
\frac{d\,J_{n}(a)}{d\,a}=-n\,J_{n+1}(a),
\end{equation*}%
\begin{equation*}
\frac{d\,^{2}J_{n}(a)}{d\,a^{2}}=n(n+1)J_{n+2}(a),
\end{equation*}%
\begin{equation*}
\frac{d^{\,3}J_{n}(a)}{da^{3}}=-n(n+1)(n+2)J_{n+3}(a).
\end{equation*}%
Here we note that a solution of the dif\thinspace ferential equation for $%
J_{1}(a)$ i.e.%
\begin{eqnarray}
&&\tfrac{1}{2}\frac{d^{\,\,3}J_{1}(a)}{d\,a^{\,3}}+2a\frac{d\,J_{1}(a)}{d\,a}%
+J_{1}(a)  \notag \\
&=&-\frac{Ai^{\,\prime }(0)^{2}}{a}-\frac{Ai(0)\,Ai^{\,\prime }(0)}{a^{2}}-%
\frac{Ai(0)^{2}}{\,a^{3}}  \label{eq21n}
\end{eqnarray}%
is suf\thinspace f\thinspace icient to obtain $\mathfrak{I}^{\,\left(
2\right) }$ since $\mathbf{J}(a)$ can be written in terms of$\ J_{1}(a)$ and
its derivatives and the relation 
\begin{equation*}
J_{1}^{\,\,\prime }(a)=-\tfrac{1}{a}Ai(0)\,Ai^{\,\prime }(0)-\tfrac{1}{2a^{2}%
}Ai(0)^{2}-J_{0}(a)+a\,J_{1}(a)+J_{3}(a).
\end{equation*}%
We get%
\begin{equation*}
\mathbf{J}(\left| a_{k}^{\prime }\right| )=-\frac{2\,Ai(0)^{2}}{5\left|
a_{k}^{\prime }\right| }+2\left| a_{k}^{\prime }\right| ^{2}[J_{1}(\left|
a_{k}^{\prime }\right| )-\frac{Ai(0)Ai^{\,\,\prime }(0)}{3\left|
a_{k}^{\prime }\right| ^{2}}-\frac{Ai^{\,\prime }(0)^{2}}{\left|
a_{k}^{\prime }\right| }]
\end{equation*}%
\begin{equation*}
+\frac{d\,J_{1}(a)}{d\,a}\mid _{a=\left| a_{k}^{\prime }\right| }+\frac{%
\left| a_{k}^{\prime }\right| }{2}\frac{d^{\,\,2}\,J_{1}(a)}{d\,a^{\,2}}\mid
_{a=\left| a_{k}^{\prime }\right| }.
\end{equation*}

\ Before proceeding to the solution of (21), we examine the behavior of $%
\mathbf{J}(\left| a_{k}^{\prime }\right| )$ for large values of $\left|
a_{k}^{\prime }\right| $. \ This is necessary given the presence of positive
powers of $\left| a_{k}^{\prime }\right| $ in the expression for $\mathbf{J}%
(\left| a_{k}^{\prime }\right| )$ shown above. \ Repeated integration by
parts of the integrals $J_{1}(a),\,J_{1}^{\,\,\prime }(a),$ and $J_{2}(a)$
give%
\begin{eqnarray*}
J_{1}(a) &\thicksim &\frac{Ai^{\,\prime }(0)^{2}}{a}+\frac{%
Ai(0)\,Ai^{\,\prime }(0)}{3a^{2}}+\frac{Ai(0)^{2}}{5a^{3}}+\cdots , \\
J_{1}^{\,\prime }(a) &\thicksim &-\frac{2Ai(0)\,Ai^{\,\prime }(0)}{3a}-\frac{%
3Ai(0)^{2}}{10a^{2}}+\frac{4Ai^{\,\prime }(0)^{2}}{7a^{3}}+\ldots , \\
J_{2}(a) &\thicksim &\frac{Ai^{\,\prime }(0)^{2}}{a^{2}}+\frac{%
2Ai(0)\,Ai^{\,\prime }(0)}{3a^{3}}+\frac{3Ai(0)^{2}}{5a^{4}}+\cdots
\end{eqnarray*}%
Using these expressions we f\thinspace ind the value of the summand $\mathbf{%
J}(\left| a_{k}^{\prime }\right| )$ as $k\rightarrow \infty $ to be%
\begin{equation*}
\mathbf{J}(\left| a_{k}^{\prime }\right| )\thicksim -\frac{3}{7\left|
a_{k}^{\prime }\right| ^{2}}\,Ai^{\,\prime }(0)^{2}-\frac{2}{3\left|
a_{k}^{\prime }\right| ^{3}}\,Ai(0)\,Ai^{\,\prime }(0)+\cdots
\end{equation*}%
and note that a cancellation of terms containing positive powers of $\left|
a_{k}^{\prime }\right| $ has taken place and the expression for the sum $%
\mathfrak{I}^{\,\left( 2\right) }$\ is seen to be f\thinspace inite for
large $k$.

The solution of the dif\thinspace ferential equation\ for $J_{1}(a)$ is then
by

\begin{eqnarray*}
\frac{J_{1}(a)}{\pi ^{2}} &=&Ai(-a)^{2}\,\left[ c_{1}-\Delta u_{1}(a)\right]
+Bi(-a)^{2}\,\left[ c_{2}-\Delta u_{2}(a)\right] \\
&&+Ai(-a)Bi(-a)\,\left[ c_{3}+2\,\Delta u_{3}(a)\right] ,
\end{eqnarray*}%
where the $c$ 's are constants of integration and 
\begin{eqnarray*}
\Delta u_{1}(a) &=&\int_{a_{0}}^{a}\left[ Ai^{\,\prime
}(0)^{2}/z+Ai(0)Ai^{\prime }(0)/z^{2}+Ai(0)^{2}/z^{3}\right]
\,Bi(-z)^{2}\,dz, \\
\Delta u_{2}(a) &=&\int_{a_{0}}^{a}\left[ Ai^{\,\prime
}(0)^{2}/z+Ai(0)Ai^{\prime }(0)/z^{2}+Ai(0)^{2}/z^{3}\right]
\,Ai(-z)^{2}\,dz, \\
\Lambda u_{3}(a) &=&\int_{a_{0}}^{a}\left[ Ai^{\,\prime
}(0)^{2}/z+Ai(0)Ai^{\prime }(0)/z^{2}+Ai(0)^{2}/z^{3}\right]
\,Ai(-z)\,Bi(-z)\,dz.
\end{eqnarray*}%
The initial conditions 
\begin{eqnarray*}
J_{1}(a_{0}) &=&J_{1}(a)|_{a=a_{0}}, \\
J_{2}(a_{0}) &=&-\frac{dJ_{1}(a)}{da}|_{a=a_{0}}\,, \\
J_{3}(a_{0}) &=&\frac{1}{2}\frac{d^{\,2}J_{1}(a)}{d\,a^{2}}|_{a=a_{0}}\,,
\end{eqnarray*}%
where $a_{0\text{ }}$is any non-zero value of $a$ are suf\thinspace
f\thinspace icient to determine the constants of integration $%
c_{1},c_{2},c_{3}$. \ The latter are given by%
\begin{eqnarray*}
c_{1}(a_{0}) &=&\,J_{1}(a_{0})[a_{0}\,Bi(-a_{0})^{2}+Bi^{\,\prime
}(-a_{0})^{2}]-\,J_{2}(a_{0})Bi(-a_{0})Bi^{\,\prime }(-a_{0}) \\
&&+\,J_{3}(a_{0})Bi(-a_{0})^{2}, \\
c_{2}(a_{0}) &=&J_{1}(a_{0})[a_{0}\,Ai(-a_{0})^{2}+Ai^{\,\prime
}(-a_{0})^{2}]\,-J_{2}(a_{0})Ai(-a_{0})Ai^{\,\prime }(-a_{0})\, \\
&&+J_{3}(a_{0})Ai(-a_{0})^{2}\,, \\
c_{3}(a_{0}) &=&-2\,J_{1}(a_{0})[a_{0}Ai(-a_{0})\,Bi(-a_{0})+Ai^{\,\prime
}(-a_{0})Bi^{\,\prime }(-a_{0})]\, \\
&&+J_{2}(a_{0})[Ai(-a_{0})Bi^{\,\prime }(-a_{0})+Ai^{\,\prime
}(-a_{0})Bi(-a_{0})]\, \\
&&-2\,J_{3}(a_{0})\,Ai(-a_{0})Bi(-a_{0})\,.
\end{eqnarray*}%
Using the relations $a_{0}=\left| a_{1}^{\,\prime }\right| ,$ $Ai^{\,\prime
}(a_{1}^{\prime })=0$ where $a_{1}^{\,\prime }$ is the f\thinspace irst root
of $Ai^{\,\prime }(a)$ and the Wronskian of $Ai(z)$ and $Bi(a),$ we have $%
Ai(a_{1}^{\prime })Bi^{\,\prime }(a_{1}^{\prime })=1/\pi $ and the
expressions for $c_{1},c_{2},c_{3}$ simplify to 
\begin{eqnarray*}
c_{1}(\left| a_{1}^{\,\prime }\right| ) &=&J_{1,3}(\left| a_{1}^{\,\prime
}\right| )\,Bi(a_{1}^{\,\prime })^{2}+[J_{1}(\left| a_{1}^{\,\prime }\right|
)\,Bi^{\text{\thinspace }\prime }(a_{1}^{\,\prime })-J_{2}(\left|
a_{1}^{\,\prime }\right| )\,Bi(a_{1}^{\,\prime })]\,Bi^{\,\prime
}(a_{1}^{\,\prime }), \\
c_{2}(\left| a_{1}^{\,\prime }\right| ) &=&J_{1,3}(\left| a_{1}^{\,\prime
}\right| )\,Ai(a_{1}^{\,\prime })^{2}, \\
c_{3}(\left| a_{1}^{\,\prime }\right| ) &=&-2\,J_{1,3}(\left|
a_{1}^{\,\prime }\right| )\,Ai(a_{1}^{\,\prime })Bi(a_{1}^{\,\prime })+\frac{%
1}{\pi }J_{2}(\left| a_{1}^{\,\prime }\right| ),
\end{eqnarray*}%
where%
\begin{equation*}
J_{1,3}(\left| a_{1}^{\,\prime }\right| )=\left| a_{1}^{\,\prime }\right|
\,J_{1}(\left| a_{1}^{\,\prime }\right| )+J_{3}(\left| a_{1}^{\,\prime
}\right| ).
\end{equation*}

\subsection{Evaluation of the $u_{i}(a)$ integrals}

The integrals $u_{i}(a)$ can be computed by expanding the functions $%
Ai(-a)^{2}$, $Bi(-a)^{2}$ and $Ai(-a)\,Bi(-a)$ in power series in $a$. \ The
three functions $Ai(z)^{2}$, $Bi(z)^{2}$ and $Ai(z)\,Bi(z)$ denoted in
common by $w(z),$ are solutions of the third-order dif\thinspace ferential
equation%
\begin{equation*}
w^{(3)}(z)+4\,z\,w^{(1)}(z)+2\,w(z)=0.
\end{equation*}%
The $n$ $th$ derivative of $w(z)$ evaluate at $z=0$ is then given by the
dif\thinspace ference equation%
\begin{equation*}
(4n+2)\,w^{\,(n)}(0)+w^{\,(n+3)}(0)=0.
\end{equation*}%
The solution of which is (where $w^{(\Box )},w^{(I)},w^{(II)}$ are
constants) 
\begin{eqnarray*}
w^{(3n)}(0) &=&w^{(\Box )}\frac{(-12)^{n}}{2\pi }\Gamma (5/6)\Gamma (n+1/6),
\\
w^{(3n+1)}(0) &=&w^{(I)}\frac{(-12)^{n}}{\sqrt{\pi }}\Gamma (n+1/2), \\
w^{(3n+2)}(0) &=&w^{(II)}\frac{(-12)^{n}\,\Gamma (n+5/6)}{\Gamma (5/6)}.
\end{eqnarray*}%
The inf\thinspace inite series expansions for the three Airy products $%
Bi(-z)^{2}$, $Ai(-z)^{2}$ and $Ai(-z)\,Bi(-z)$\ are then%
\begin{eqnarray*}
Bi(-z)^{2} &=&w_{1}^{(\Box
)}\,S_{1}(z)+w_{1}^{(I)}\,S_{2}(z)+w_{1}^{(II)}\,S_{3}(z), \\
Ai(-z)^{2} &=&w_{2}^{(\Box
)}\,S_{1}(z)+w_{2}^{(I)}\,S_{2}(z)+w_{2}^{(II)}\,S_{3}(z), \\
Ai(-z)\,Bi(-z) &=&w_{3}^{(\Box
)}\,S_{1}(z)+w_{3}^{(I)}\,S_{2}(z)+w_{3}^{(II)}\,S_{3}(z),
\end{eqnarray*}%
where%
\begin{eqnarray*}
S_{1}(z) &=&\frac{\Gamma (5/6)}{2\pi }\sum_{n=0}^{\infty }\frac{%
12^{n}\,\Gamma (n+1/6)}{(3n)!}(-z)^{3n}, \\
S_{2}(z) &=&-\frac{z}{\sqrt{\pi }}\sum_{n=0}^{\infty }\frac{12^{n}\,\Gamma
(n+1/2)}{(3n+1)!}(-z)^{3n}, \\
S_{3}(z) &=&\frac{z^{2}}{\Gamma (5/6)}\sum_{n=0}^{\infty }\frac{%
12^{n}\,\Gamma (n+5/6)}{(3n+2)!}(-z)^{3n},
\end{eqnarray*}%
and $w_{i}^{(\Box )},w_{i}^{(I)},\,w_{i}^{(II)}$ are given in the Table 3
below.

\begin{center}
\begin{table}[H] \centering%
\caption{The constants w\label{four}}%
\end{table}%

\renewcommand{\arraystretch}{2.0}%
\begin{equation*}
\begin{tabular}{|c|c|c|c|c|}
\hline
$i$ & $\,w_{i}(z)$ & $\,w_{i}^{(\Box )}$ & $\,w_{i}^{(I)}$ & $w_{i}^{(II)}$
\\ \hline
$1$ & $Bi(z)^{2}$ & $3Ai^{\,}(0)^{2}$ & $-\,6\,Ai(0)\,Ai^{\,\prime }(0)$ & $%
6\,Ai^{\,\prime }(0)^{2}$ \\ \hline
$2$ & $Ai(z)^{2}$ & $Ai^{\,}(0)^{2}$ & $2\,Ai(0)\,Ai^{\,\prime }(0)$ & $%
2\,Ai^{\,\prime }(0)^{2}$ \\ \hline
$3$ & $Ai(z)Bi(z)$ & $\sqrt{3}Ai(0)^{2}$ & $0$ & $-2\sqrt{3}Ai^{\,\prime
}(0)^{2}$ \\ \hline
\end{tabular}%
\end{equation*}
\end{center}

Integrating the terms appearing in the series representation of $%
w_{i}(-z)/z^{\,n}$ and summing the resulting terms produces closed form
expression involving hypergeometric functions e.g.%
\begin{eqnarray*}
\int \frac{w_{i}(-z)}{z}dz &=&w_{i}^{(\Box )}\left[ \ln (z)-\frac{z^{3}}{9}%
\,_{3}F_{4}\tbinom{1,1,7/6\,\,;\;\frac{4}{9}z^{3}}{4/3,5/3,2,2}\right] \\
&&+z\,\text{\/}\,w_{i}^{(I)}\,\,\,_{2}F_{3}\tbinom{1/3,1/2\;\,;\;\,-\frac{4}{%
9}z^{3}}{2/3,4/3,4/3} \\
&&+z^{2}\frac{w_{i}^{(II)}\,}{4}\,_{2}F_{3}\tbinom{2/3,5/6\,\text{\/}\;;\;-%
\frac{4}{9}z^{3}}{4/3,5/3,5/3},
\end{eqnarray*}%
\begin{eqnarray*}
\int \frac{w_{i}(-z)}{z^{2}}dz &=&-\frac{w_{i}^{(\Box )}}{z}\,\,_{2}F_{3}%
\tbinom{-1/3,1/6\;\,\text{\/};\;-\frac{4}{9}z^{3}}{1/3,2/3,2/3} \\
&&+w_{i}^{(I)}\left[ \ln (z)-\frac{z^{3}}{12}\,_{3}F_{4}\tbinom{1,1,3/2\;;\;-%
\frac{4}{9}z^{3}}{2,2,5/3,7/3}\right] \\
&&+\,z\,\frac{w_{i}^{(II)}}{2}\,_{2}F_{3}\tbinom{1/3,5/6\;;\;-\frac{4}{9}%
z^{3}}{4/3,4/3,5/3},
\end{eqnarray*}%
\begin{eqnarray*}
\int \frac{w_{i}(-z)}{z^{3}}dz &=&-\frac{w_{i}^{(\Box )}}{2z^{2}}\,_{2}F_{3}%
\tbinom{-2/3,1/6\;;\;-\frac{4}{9}z^{3}}{1/3,1/3,2/3}-\frac{w_{i}^{(I)}}{z}%
\,\,_{2}F_{3}\tbinom{-1/3,1/2\;;\;-\frac{4}{9}z^{3}}{2/3,2/3,4/3} \\
&&+\frac{w_{i}^{(II)}}{2}\left[ \ln (z)-\frac{z^{3}}{18}\,\,_{3}F_{4}\tbinom{%
1,1,11/6\;;\;-\frac{4}{9}z^{3}}{2,2,7/3,8/3}\right] .
\end{eqnarray*}%
Alternately, Maple gives the same results. \ The resulting explicit
expressions for the $u_{i}(z)$ are%
\begin{eqnarray*}
u_{1}(z) &=&\frac{\pi ^{2}}{2\,\cdot 3^{8/3}\Gamma (2/3)^{4}\,z^{2}}%
\,_{2}F_{3}\tbinom{-2/3,1/6\;;\;-\frac{4}{9}z^{3}}{1/3,1/3,2/3} \\
&&-\frac{3^{2/3}\Gamma (2/3)^{4}\,z^{2}\,}{32\pi ^{2}}\,_{2}F_{3}\tbinom{%
2/3,5/6\,\text{\/}\;;\;-\frac{4}{9}z^{3}}{4/3,5/3,5/3} \\
&&+\frac{\pi }{3^{11/6}\Gamma (2/3)^{2}\,\,z}\left[ \,_{2}F_{3}\tbinom{%
-1/3,1/2\;;\;-\frac{4}{9}z^{3}}{2/3,2/3,4/3}-\frac{1}{2}\,_{2}F_{3}\tbinom{%
-1/3,1/6\;\,\text{\/};\;-\frac{4}{9}z^{3}}{1/3,2/3,2/3}\right] \\
&&+z\frac{\Gamma (2/3)^{2}}{8\pi \,3^{1/6}}\,\left[ \,_{2}F_{3}\tbinom{%
1/3,5/6\;;\;-\frac{4}{9}z^{3}}{4/3,4/3,5/3}-2\,\,\,_{2}F_{3}\tbinom{%
1/3,1/2\;\,;\;\,-\frac{4}{9}z^{3}}{2/3,4/3,4/3}\right] \\
&&+\frac{z^{3}}{216}\left[ 
\begin{array}{c}
\,_{3}F_{4}\tbinom{1,1,11/6\;;\;-\frac{4}{9}z^{3}}{2,2,7/3,8/3}+2\,_{3}F_{4}%
\tbinom{1,1,7/6\,\,;\;\frac{4}{9}z^{3}}{4/3,5/3,2,2} \\ 
-3\,_{3}F_{4}\tbinom{1,1,3/2\;;\;-\frac{4}{9}z^{3}}{2,2,5/3,7/3}%
\end{array}%
\right] , \\
u_{2}(z) &=&-\ln (z)+\frac{\pi ^{2}}{2\cdot 3^{5/3}\Gamma (2/3)^{4}\,z^{2}}%
\,_{2}F_{3}\tbinom{-2/3,1/6\;;\;-\frac{4}{9}z^{3}}{1/3,1/3,2/3} \\
&&-\frac{3^{5/3}\Gamma (2/3)^{4}}{32\pi ^{2}}z^{2}\,\,_{2}F_{3}\tbinom{%
2/3,5/6\,\text{\/}\;;\;-\frac{4}{9}z^{3}}{4/3,5/3,5/3} \\
&&-\frac{\pi }{\,3^{5/6}\Gamma (2/3)^{2}\,z}\left[ \,_{2}F_{3}\tbinom{%
-1/3,1/2\;;\;-\frac{4}{9}z^{3}}{2/3,2/3,4/3}+\frac{1}{2}\,_{2}F_{3}\tbinom{%
-1/3,1/6\;\,\text{\/};\;-\frac{4}{9}z^{3}}{1/3,2/3,2/3}\right] \\
&&+\frac{\,3^{5/6}\Gamma (2/3)^{2}\,z}{8\pi }\left[ \,\,_{2}F_{3}\tbinom{%
1/3,5/6\;;\;-\frac{4}{9}z^{3}}{4/3,4/3,5/3}+2\,_{2}F_{3}\tbinom{%
1/3,1/2\;\,;\;\,-\frac{4}{9}z^{3}}{2/3,4/3,4/3}\right] \\
&&+\frac{z^{3}}{72}\left[ 
\begin{array}{c}
\,_{3}F_{4}\tbinom{1,1,11/6\;;\;-\frac{4}{9}z^{3}}{2,2,7/3,8/3}+2\,_{3}F_{4}%
\tbinom{1,1,7/6\,\,;\;\frac{4}{9}z^{3}}{4/3,5/3,2,2} \\ 
+3\,_{3}F_{4}\tbinom{1,1,3/2\;;\;-\frac{4}{9}z^{3}}{2,2,5/3,7/3}%
\end{array}%
\right] , \\
u_{3}(z) &=&-\frac{\pi ^{2}}{3^{13/6}\Gamma (2/3)^{4}\,z^{2}}\,\,_{2}F_{3}%
\tbinom{-2/3,1/6\;;\;-\frac{4}{9}z^{3}}{1/3,1/3,2/3}- \\
&&\frac{3^{7/6}\Gamma (2/3)^{4}}{16\pi ^{2}}\,z^{2}\,_{2}F_{3}\tbinom{%
2/3,5/6\,\text{\/}\;;\;-\frac{4}{9}z^{3}}{4/3,5/3,5/3} \\
&&+\frac{\pi }{3^{4/3}\Gamma (2/3)^{2}\,z}\,_{2}F_{3}\tbinom{-1/3,1/6\;\,%
\text{\/};\;-\frac{4}{9}z^{3}}{1/3,2/3,2/3}+ \\
&&\frac{3^{1/3}\,\Gamma (2/3)^{2}}{4\pi }\,z\,\,_{2}F_{3}\tbinom{%
1/3,5/6\;;\;-\frac{4}{9}z^{3}}{4/3,4/3,5/3} \\
&&+\frac{\sqrt{3}z^{3}}{108}\left[ \,_{3}F_{4}\tbinom{1,1,11/6\;;\;-\frac{4}{%
9}z^{3}}{2,2,7/3,8/3}-2\,_{3}F_{4}\tbinom{1,1,7/6\,\,;\;\frac{4}{9}z^{3}}{%
4/3,5/3,2,2}\right] .
\end{eqnarray*}

\subsection{The $\mathbf{J(}a\mathbf{)}$ terms}

\bigskip The f\thinspace irst and second derivatives of $J_{1}\ (a)$ are in
given by%
\begin{eqnarray*}
\frac{1}{\pi ^{2}}\frac{d\,\,J_{1}\ (a)}{d\,a} &=&-\,2\,Ai(-a)Ai^{\,\prime
}(-a)[c_{1}-\Delta u_{1}(a)] \\
&&-\,2\,Bi(-a)Bi^{\,\prime }(-a)[c_{2}-\Delta u_{2}(a)] \\
&&-[Ai^{\,\prime }(-a)Bi(-a)+Ai(-a)Bi^{\,\prime }(-a)][c_{3}+2\Delta
\,u_{3}(a)],
\end{eqnarray*}%
\begin{eqnarray*}
\frac{1}{\pi ^{2}}\frac{d\,^{2}\,J_{1}\ (a)}{d\,a^{2}} &=&2\,[Ai^{\,\prime
}(-a)^{2}-a\,Ai(-a)^{2}][c_{1}-\Delta u_{1}(a)] \\
&&+2\,[Bi^{\,\prime }(-a)^{2}-a\,Bi(-a)^{2}][c_{2}-\Delta u_{2}(a)] \\
&&-2\,[a\,Ai(-a)Bi(-a)-Ai^{\,\prime }(-a)Bi^{\,\prime }(-a)][c_{3}+2\,\Delta
u_{3}(a)].
\end{eqnarray*}%
Using these expressions, the terms $\mathbf{J}(\left| a_{k}^{\prime }\right|
)$ are 
\begin{equation*}
\mathbf{J}(\left| a_{k}^{\prime }\right| )=-\frac{2\,Ai(0)^{2}}{5\left|
a_{k}^{\prime }\right| }-\frac{2}{3}Ai(0)Ai^{\,\,\prime }(0)-2\left|
a_{k}^{\prime }\right| Ai^{\,\prime }(0)^{2}+\mathsf{j}(a)|_{a=\left|
a_{k}^{\prime }\right| },
\end{equation*}%
where%
\begin{equation*}
\mathsf{j}(a)=\frac{1}{\pi ^{2}}[2a^{2}\,J_{1}(a)+\frac{d\,\,J_{1}\ (a)}{d\,a%
}+\frac{a}{2}\frac{d\,^{2}\,J_{1}\ (a)}{d\,a^{2}}],
\end{equation*}%
and $\mathsf{j}(a)|_{a=\left| a_{k}^{\prime }\right| }$ hereafter denoted by 
$\mathsf{j}(\left| a_{k}^{\prime }\right| )$ is explicitly given by 
\begin{eqnarray*}
&&\mathsf{j}(\left| a_{k}^{\prime }\right| )=[\left| a_{k}^{\prime }\right|
^{2}Ai(a_{k}^{\prime })^{2}-2Ai(a_{k}^{\prime })Ai^{\,\prime }(a_{k}^{\prime
})+\left| a_{k}^{\prime }\right| Ai^{\,\prime }(a_{k}^{\prime })^{2}] \\
&&\cdot \lbrack c_{1}-\Delta u_{1}(\left| a_{k}^{\prime }\right| )] \\
&&+[\left| a_{k}^{\prime }\right| ^{2}Bi(a_{k}^{\prime
})^{2}-2Bi(a_{k}^{\prime })Bi^{\,\prime }(a_{k}^{\prime })+\left|
a_{k}^{\prime }\right| Bi^{\,\prime }(a_{k}^{\prime })^{2}] \\
&&\cdot \lbrack c_{2}-\Delta u_{2}(\left| a_{k}^{\prime }\right| )] \\
&&+[c_{3}+2\,\Delta u_{3}(\left| a_{k}^{\prime }\right| )]\cdot \lbrack
\,\left| a_{k}^{\prime }\right| ^{2}Ai(a_{k}^{\prime })Ai^{\,\prime
}(a_{k}^{\prime })+\left| a_{k}^{\prime }\right| Ai^{\,\prime
}(a_{k}^{\prime })Bi^{\,\prime }(a_{k}^{\prime }) \\
&&-Ai(a_{k}^{\prime })Bi^{\,\prime }(a_{k}^{\prime })-Ai^{\,\prime
}(a_{k}^{\prime })Bi(a_{k}^{\prime })],
\end{eqnarray*}%
which simplif\thinspace ies to%
\begin{eqnarray*}
\mathsf{j}(\left| a_{k}^{\prime }\right| ) &=&[\,\left| a_{k}^{\prime
}\right| ^{2}Ai(a_{k}^{\prime })^{2}][c_{1}-\Delta u_{1}(\left|
a_{k}^{\prime }\right| )] \\
&&+[\,\left| a_{k}^{\prime }\right| ^{2}Bi(a_{k}^{\prime
})^{2}-2Bi(a_{k}^{\prime })Bi\,^{\prime }(a_{k}^{\prime })+\left|
a_{k}^{\prime }\right| \,Bi^{\,\prime }(a_{k}^{\prime })^{2}] \\
&&\cdot \lbrack c_{2}-\Delta u_{2}(\left| a_{k}^{\prime }\right| )] \\
&&-[c_{3}+2\,\Delta u_{3}(\left| a_{k}^{\prime }\right| )]\,Ai(a_{k}^{\prime
})Bi^{\,\prime }(a_{k}^{\prime }).
\end{eqnarray*}%
Upon insertion of the expressions for $\Delta u_{i}(\left| a_{k}^{\prime
}\right| )$ and collecting terms we get the complicated form for $\mathsf{j}%
(\left| a_{k}^{\prime }\right| )$\ 
\begin{eqnarray*}
\mathsf{j}(\left| a_{k}^{\prime }\right| ) &=&c_{1}\,[\,\left| a_{k}^{\prime
}\right| ^{2}Ai(a_{k}^{\prime })^{2}-2Ai(a_{k}^{\prime })Ai^{\,\prime
}(a_{k}^{\prime })+\left| a_{k}^{\prime }\right| Ai^{\,\prime
}(a_{k}^{\prime })^{2}] \\
&&+[c_{2}+\ln (\left| a_{k}^{\prime }/a_{1}^{\prime }\right| )\,] \\
&&\cdot \lbrack \left| a_{k}^{\prime }\right| ^{2}Bi(a_{k}^{\prime
})^{2}-2Bi(a_{k}^{\prime })Bi^{\,\prime }(a_{k}^{\prime })+\left|
a_{k}^{\prime }\right| Bi^{\,\prime }(-a)^{2}] \\
&&+c_{3}\,[\,\left| a_{k}^{\prime }\right| ^{2}Ai(a_{k}^{\prime
})Ai^{\,\prime }(a_{k}^{\prime })+\left| a_{k}^{\prime }\right| Ai^{\,\prime
}(a_{k}^{\prime })Bi^{\,\prime }(a_{k}^{\prime }) \\
&&-Ai(a_{k}^{\prime })Bi^{\,\prime }(a_{k}^{\prime })-Ai^{\,\prime
}(a_{k}^{\prime })Bi(a_{k}^{\prime })] \\
&&-d_{1}(\left| a_{k}^{\prime }\right| )\left[ 
\begin{array}{c}
-\frac{\pi ^{2}}{2\cdot 3^{2/3}\Gamma (2/3)^{4}}\Delta \,\{_{2}F_{3}\tbinom{%
-2/3,1/6\;;\;-\frac{4}{9}a^{3}}{1/3,1/3,2/3}/a^{2}\} \\ 
+\frac{\pi }{2\cdot 3^{11/6}\Gamma (2/3)^{2}\,}\Delta \,\{_{2}F_{3}\tbinom{%
-1/3,1/6\;\,\text{\/};\;-\frac{4}{9}a^{3}}{1/3,2/3,2/3}\,/a\} \\ 
-\frac{1}{108\,}\Delta \,\{a^{3}\,_{3}F_{4}\tbinom{1,1,7/6\,\,;\;\frac{4}{9}%
a^{3}}{4/3,5/3,2,2}\}%
\end{array}%
\right] \\
&&-d_{2}(\left| a_{k}^{\prime }\right| )\left[ 
\begin{array}{c}
\frac{3^{2/3}\Gamma (2/3)^{4}}{32\pi ^{2}}\Delta \{\,a^{2}\,_{2}F_{3}\tbinom{%
2/3,5/6\,\text{\/}\;;\;-\frac{4}{9}a^{3}}{4/3,5/3,5/3}\} \\ 
-\frac{\Gamma (2/3)^{2}\,}{8\pi \,3^{1/6}}\Delta \,\{a\,\,_{2}F_{3}\tbinom{%
1/3,5/6\;;\;-\frac{4}{9}a^{3}}{4/3,4/3,5/3}\} \\ 
-\frac{1}{216}\Delta \,\{a^{3}\,_{3}F_{4}\tbinom{1,1,11/6\;;\;-\frac{4}{9}%
a^{3}}{2,2,7/3,8/3}\}%
\end{array}%
\right] \\
&&+d_{3}(\left| a_{k}^{\prime }\right| )\left[ 
\begin{array}{c}
\frac{\pi }{3^{11/6}\Gamma (2/3)^{2}\,}\Delta \{_{2}F_{3}\tbinom{%
-1/3,1/2\;;\;-\frac{4}{9}a^{3}}{2/3,2/3,4/3}/a\} \\ 
-\frac{\Gamma (2/3)^{2}}{4\pi \,3^{1/6}}\,\Delta \{\,a\,_{2}F_{3}\tbinom{%
1/3,1/2\;\,;\;\,-\frac{4}{9}a^{3}}{2/3,4/3,4/3}\} \\ 
-\frac{1}{72}\Delta \{\,a^{3}\,_{3}F_{4}\tbinom{1,1,3/2\;;\;-\frac{4}{9}a^{3}%
}{2,2,5/3,7/3}\}%
\end{array}%
\right] ,
\end{eqnarray*}%
where%
\begin{eqnarray*}
d_{1}(a) &=&a^{2}[Ai(-a)^{2}+4\sqrt{3}Ai(-a)Bi(-a)+Bi(-a)^{2}] \\
&&+[aAi^{\,\prime }(-a)^{2}+4\sqrt{3}Ai^{\,\prime }(-a)Bi^{\,\prime
}(-a)+Bi^{\,\prime }(-a)^{2}] \\
&&-4\sqrt{3}\left[ Ai(-a)Bi^{\,\prime }(-a)+Ai^{\,\prime }(-a)Bi(-a)\right]
\\
&&-2\left[ Ai(-a)Ai^{\,\prime }(-a)+3Bi(-a)Bi^{\,\prime }(-a)\right] , \\
&& \\
d_{2}(a) &=&a^{2}[Ai(-a)^{2}-4\sqrt{3}Ai(-a)Bi(-a)+Bi(-a)^{2}] \\
&&+a[Ai^{\,\prime }(-a)^{2}-4\sqrt{3}Ai^{\,\prime }(-a)Bi^{\,\prime
}(-a)+Bi^{\,\prime }(-a)^{2}] \\
&&+4\sqrt{3}\left[ Ai(-a)Bi^{\,\prime }(-a)+Ai^{\,\prime }(-a)Bi(-a)\right]
\\
&&-2\left[ Ai(-a)Ai^{\,\prime }(-a)+3Bi(-a)Bi^{\,\prime }(-a)\right] , \\
&& \\
d_{3}(a) &=&a^{2}[Ai(-a)^{2}-3Bi(-a)^{2}]+a[Ai^{\,\prime
}(-a)^{2}-3Bi^{\,\prime }(-a)^{2}] \\
&&-2[Ai(-a)Ai^{\,\prime }(-a)-3Bi(-a)Bi^{\,\prime }(-a)].
\end{eqnarray*}%
In the case where $a=\left| a_{k}^{\prime }\right| $ we have $%
Ai(a_{k}^{\prime })Bi^{\prime }(a_{k}^{\prime })=1/\pi $ and the expressions
for $\mathsf{j}(\left| a_{k}^{\prime }\right| )$ and $d_{i}(a)$ become%
\begin{eqnarray*}
\mathsf{j}(\left| a_{k}^{\prime }\right| ) &=&c_{1}\,[\,\left| a_{k}^{\prime
}\right| ^{2}Ai(a_{k}^{\prime })^{2}]-c_{3}/\pi \\
&&+[c_{2}+\ln (a_{k}^{\prime }/a_{1}^{\prime })]\,[\,\left| a_{k}^{\prime
}\right| ^{2}Bi(a_{k}^{\prime })^{2}-2Bi(a_{k}^{\prime })Bi^{\,\prime
}(a_{k}^{\prime })+\left| a_{k}^{\prime }\right| Bi^{\,\prime
}(a_{k}^{\prime })^{2}] \\
&&-d_{1}(\left| a_{k}^{\prime }\right| )\,\Delta \widetilde{H}%
_{1}-d_{2}(\left| a_{k}^{\prime }\right| )\,\Delta \widetilde{H}%
_{2}+d_{3}(\left| a_{k}^{\prime }\right| )\,\Delta \widetilde{H}_{3},
\end{eqnarray*}%
and where%
\begin{eqnarray*}
\widetilde{H}_{1}(a) &=&-\tfrac{\pi ^{2}}{2\cdot 3^{2/3}\Gamma (2/3)^{4}a^{2}%
}\,\text{ }_{2}F_{3}\tbinom{-2/3,1/6\;;\;-\frac{4}{9}a^{3}}{1/3,1/3,2/3}-%
\tfrac{1}{108\,}\,a^{3}\,_{3}F_{4}\tbinom{1,1,7/6\,\,;\;\frac{4}{9}a^{3}}{%
4/3,5/3,2,2} \\
&&+\tfrac{\pi }{2\cdot 3^{11/6}\Gamma (2/3)^{2}\,a}\,_{2}F_{3}\tbinom{%
-1/3,1/6\;\,\text{\/};\;-\frac{4}{9}a^{3}}{1/3,2/3,2/3}, \\
\widetilde{H}_{2}(a) &=&\tfrac{3^{2/3}\Gamma (2/3)^{4}}{32\pi ^{2}}%
\,a^{2}\,_{2}F_{3}\tbinom{2/3,5/6\,\text{\/}\;;\;-\frac{4}{9}a^{3}}{%
4/3,5/3,5/3}-\tfrac{\Gamma (2/3)^{2}\,}{8\pi \,3^{1/6}}\,a\,\,_{2}F_{3}%
\tbinom{1/3,5/6\;;\;-\frac{4}{9}a^{3}}{4/3,4/3,5/3} \\
&&-\tfrac{a^{3}}{216}\,\,_{3}F_{4}\tbinom{1,1,11/6\;;\;-\frac{4}{9}a^{3}}{%
2,2,7/3,8/3}, \\
\widetilde{H}_{3}(a) &=&\tfrac{\pi }{3^{11/6}\Gamma (2/3)^{2}\,a}%
\,\,_{2}F_{3}\tbinom{-1/3,1/2\;;\;-\frac{4}{9}a^{3}}{2/3,2/3,4/3}-\tfrac{%
\Gamma (2/3)^{2}}{4\pi \,3^{1/6}}\,\,a\,\,_{2}F_{3}\tbinom{1/3,1/2\;\,;\;\,-%
\frac{4}{9}a^{3}}{2/3,4/3,4/3} \\
&&-\tfrac{a^{3}}{72}\,_{3}F_{4}\tbinom{1,1,3/2\;;\;-\frac{4}{9}a^{3}}{%
2,2,5/3,7/3},
\end{eqnarray*}%
and

\begin{eqnarray*}
d_{1}(\left| a_{k}^{\prime }\right| ) &=&\left| a_{k}^{\prime }\right|
^{2}[Ai(a_{k}^{\prime })^{2}+4\sqrt{3}Ai(a_{k}^{\prime })Bi(a_{k}^{\prime
})+Bi(a_{k}^{\prime })^{2}] \\
&&+\left| a_{k}^{\prime }\right| Bi^{\,\prime }(a_{k}^{\prime
})^{2}-6Bi(a_{k}^{\prime })Bi^{\,\prime }(a_{k}^{\prime })-4\sqrt{3}/\pi , \\
d_{2}(\left| a_{k}^{\prime }\right| ) &=&\left| a_{k}^{\prime }\right|
^{2}[Ai(a_{k}^{\prime })^{2}-4\sqrt{3}Ai(a_{k}^{\prime })Bi(a_{k}^{\prime
})+Bi(a_{k}^{\prime })^{2}] \\
&&+\left| a_{k}^{\prime }\right| Bi^{\,\prime }(a_{k}^{\prime
})^{2}-6Bi(a_{k}^{\prime })Bi^{\,\prime }(a_{k}^{\prime })+4\sqrt{3}/\pi , \\
d_{3}(\left| a_{k}^{\prime }\right| ) &=&\left| a_{k}^{\prime }\right|
^{2}[Ai(a_{k}^{\prime })^{2}-3Bi(a_{k}^{\prime })^{2}]-3\left| a_{k}^{\prime
}\right| Bi^{\,\prime }(a_{k}^{\prime })^{2}+6Bi(a_{k}^{\prime
})Bi^{\,\prime }(a_{k}^{\prime }).
\end{eqnarray*}

As in the case of $\mathfrak{I}^{\left( 1\right) }$ where the rate of
convergence of the sum over $\mathbf{I}_{1}$ could be accelerated if the sum
containing $\left| a_{k}^{\prime }\right| $ with large $k$ was subtracted, a
similar procedure can be performed in the case of\ $\mathfrak{I}^{\,\left(
2\right) }$. \ We write for $\mathfrak{I}^{\left( 2\right) }$ \ \textbf{\ }%
\begin{equation*}
\mathfrak{I}^{\left( 2\right) }=\tfrac{1}{3Ai^{\,\prime }(0)^{2}}%
\sum_{k=1}^{N}[\mathbf{J}(\left| a_{k}^{\prime }\right| )-\mathbf{J}(\left|
a_{k}^{\prime }\right| )_{\left| a_{i}^{\prime }\right| >x}]+\tfrac{1}{%
3Ai^{\,\prime }(0)^{2}}\sum_{i=N+1}^{\infty }\mathbf{J}(\left| a_{k}^{\prime
}\right| )_{\left| a_{i}^{\prime }\right| >x},
\end{equation*}%
and def\thinspace ining%
\begin{eqnarray*}
\mathbf{J}(a)_{a\,\geq \,\,x} &=&\frac{1}{a^{2}}\int_{0}^{a}\frac{x}{(1+x/a)}%
[2\,Ai(x)Ai^{\,\prime }(x)+xAi^{\,\prime }(x)^{2}-x^{2}Ai(x)^{2}], \\
&& \\
\mathbf{J}(a)_{a\,\,\geq \,\,x} &\simeq &\frac{1}{a^{2}}\int_{0}^{\infty }%
\frac{x}{(1+x/a)}[2\,Ai(x)Ai^{\,\prime }(x)+xAi^{\,\prime
}(x)^{2}-x^{2}Ai(x)^{2}],
\end{eqnarray*}%
where in the last expression, terms of order $\exp (-4a^{3/2}/3)$ have been
omitted. \ Expanding the denominator in the expression above in a power
series in $x$ followed by integration (cf. Appendix III) gives the
approximate expression for $\mathfrak{I}^{\left( 2\right) }$%
\begin{eqnarray*}
\mathfrak{I}^{\left( 2\right) } &\simeq &\tfrac{1}{3Ai^{\,\prime }(0)^{2}}%
\sum_{k=1}^{N}\mathbf{J}(\left| a_{k}^{\,\prime }\right| ) \\
&&-\tfrac{2}{12^{13/6}\sqrt{\pi }Ai^{\,\prime }(0)^{2}}\sum_{k=0}^{n}\tfrac{%
\left( -1\right) ^{k}}{12^{\,k/3}}\tfrac{(k+4)(k+1)!}{\Gamma (k/3+13/6)}%
\left\{ \mathcal{Z}_{k+2}-\mathcal{Z}_{k+2}(N)\right\} .
\end{eqnarray*}%
In the case where $N=10$ and $n=6$ the expression above produces $\mathfrak{I%
}^{\left( 2\right) }=-\ 0.2636317121$ a value which is accurate to 8 decimal
places. \ 

\bigskip

\subsection{Values of $J_{n}(a)$ for small $a$}

The solution of the third-order dif\thinspace f\thinspace erential equation\
(21) for $J_{1}(a)$ requires values of $J_{1}(a_{0}),$ $J_{2}(a_{0}),$ and $%
J_{3}(a_{0}).$ \ If analytic expressions for those quantities are also
sought, then the integrals $J_{n}(a)$ involving $Ai(x)^{2}$ can be obtained
in a way similar to that used for the integrals $I_{n}(a)$ which contained $%
Ai(x)$ and were treated above. \ That is to say we write 
\begin{equation}
J_{n}(a)=\int_{a}^{\infty }\frac{Ai(x-a)^{2}}{x^{\,n}}dx,  \label{eq22}
\end{equation}%
and expand $Ai(x-a)^{2}$ in a power series in $a$. We get%
\begin{equation*}
Ai(x-a)^{2}=\sum_{j=0}^{\infty }\frac{(-a)^{\,\,j}}{j!}\frac{%
d^{\,\,\,j}\,Ai(x)^{2}}{d\,x^{\,\,\,j}}
\end{equation*}%
\ \ To begin with, we compute the $j$ $th$ derivatives of $Ai(x)^{2}.$ \ It
is easy to see that these derivatives are given by expressions with the
following form%
\begin{equation*}
\frac{d^{\;\,j}\,Ai(x)^{2}}{d\,x^{\;j\,}}=\mathfrak{P}_{j}(x)\,Ai(x)^{2}+%
\mathfrak{Q}_{\,j}(x)\,Ai^{\,\prime }(x)^{2}+\mathfrak{R\,}%
_{j}(x)\,Ai(x)\,Ai^{\,\prime }(x),
\end{equation*}%
where the quantities $\mathfrak{P}_{j}(x),$ $\mathfrak{Q}_{\,\,j}(x),$%
\bigskip $\,\mathfrak{R}\,_{j}(x)$ are polynomials. \ Using these relations
we have 
\begin{equation*}
Ai(x-a)^{2}=\sum_{j=0}^{\infty }\frac{(-a)^{\,j}}{j!}\left[ \mathfrak{P}%
_{j}(x)\,Ai(x)^{2}+\mathfrak{Q}_{\,j}(x)\,Ai^{\,\prime }(x)^{2}+\mathfrak{R\,%
}\,_{j}(x)\,Ai(x)\,Ai^{\,\prime }(x)\right] .
\end{equation*}%
Rearranging the latter expression we get 
\begin{eqnarray*}
Ai(x-a)^{2} &=&Ai(x)^{2}\sum_{j=0}^{\infty }\frac{(-a)^{\,j}}{j!}\,\mathfrak{%
P\,}_{j}(x)\,+Ai^{\,\prime }(x)^{2}\sum_{j=0}^{\infty }\frac{(-a)^{\,j}}{j!}%
\,\mathfrak{Q\,}_{\,j}(x)\, \\
&&+Ai(x)\,Ai^{\,\prime }(x)\sum_{j=0}^{\infty }\frac{(-a)^{\,j}}{j!}\,%
\mathfrak{R}\,_{j}(x)\,.
\end{eqnarray*}%
Expansions of $Ai(x-a)^{2}$ in a Taylor series which includes powers of $a$
up to $a^{15}$, allows values of $J_{1}(\left| a_{1}^{\prime }\right| ),\
J_{2}(\left| a_{1}^{\prime }\right| )$ and $J_{3}(\left| a_{1}^{\prime
}\right| )\ $to be obtained which are accurate to\textbf{\ }8 decimal places
with values of 0.04826441, 0.03654795 and 0.02879280 respectively. \ As was
seen above this method is not easily expressed in general terms and a more
systematic way of obtaining an analytical expression for the $J_{n}(a)$
integrals is given below.\medskip

The quantities in the brackets above have the form of generating functions
and are here denoted by \ $\Xi (-a,x),\,\Lambda (-a,x)$, and $\varrho (-a,x)$
for the polynomials $\mathfrak{P}_{j}(x),$ $\mathfrak{Q}_{\,\,j}(x),$\ and $%
\mathfrak{R}\,_{j}(x)$ respectively. \ Written more compactly we have 
\begin{equation}
J_{n}(a)=\int_{a}^{\infty }\left[ \Xi (-a,x)\,Ai(x)^{2}+\Lambda
(-a,x)\,Ai^{\,\prime }(x)^{2}+\varrho \mathfrak{\,}(-a,x)\,Ai(x)\,Ai^{\,%
\prime }(x)\right] \frac{dx}{x^{\,n}}.  \label{eq23}
\end{equation}

From the expression corresponding to the $j+1$ derivative of $Ai(x)^{2}$ we
get the dif\thinspace ferential recurrence equations for the $\mathfrak{P}%
_{j}(x),\mathfrak{Q}_{\,j}(x)$ and $\mathfrak{R\,}_{j}(x)$ polynomials i.e. 
\begin{eqnarray*}
\mathfrak{P}_{j+1}(x) &=&\frac{d\,\mathfrak{P}_{j}(x)}{d\,x}+x\,\mathfrak{R\,%
}_{j}(x), \\
\mathfrak{Q}_{\,j+1}(x) &=&\frac{d\,\mathfrak{Q}_{\,j}(x)}{d\,x}+\mathfrak{%
R\,}\,_{j}(x), \\
\mathfrak{R\,}\,_{j+1}(x) &=&\frac{d\,\mathfrak{R\,}\,_{j}(x)}{d\,x}+2%
\mathfrak{P}\,_{j}(x)+2\,x\,\mathfrak{Q}_{\,j}(x),
\end{eqnarray*}%
with initial values given by 
\begin{eqnarray*}
\mathfrak{P}_{0}(x) &=&1, \\
\mathfrak{Q}_{0}(x) &=&0, \\
\mathfrak{R\,}_{0}(x) &=&0.
\end{eqnarray*}%
The f\thinspace irst few of these polynomials are given in the Table 4 below.

\bigskip

\begin{center}
\begin{table}[H] \centering%
\caption{The polynomials P, Q and R\label{five}}%
\end{table}%

\begin{tabular}{|c|c|c|c|}
\hline
$i$ & $\mathfrak{P}_{i}$ & $\mathfrak{Q}_{i}$ & $\mathfrak{R\,}_{i}$ \\ 
\hline
$0$ & $1$ & $0$ & $0$ \\ \hline
$1$ & $0$ & $0$ & $2$ \\ \hline
$2$ & $2\,x$ & $2$ & $0$ \\ \hline
$3$ & $2$ & $0$ & $8\,x$ \\ \hline
$4$ & $8\,x^{2}$ & $8\,x$ & $12$ \\ \hline
$5$ & $28\,x$ & $20$ & $32\,x^{2}$ \\ \hline
$6$ & $28+32\,x^{3}$ & $32\,x^{2}$ & $160\,x$ \\ \hline
$7$ & $256\,x^{2}$ & $224\,x$ & $216+128\,x^{3}$ \\ \hline
$8$ & $728\,x+128\,x^{4}$ & $440+128\,x^{3}$ & $1344\,x^{2}$ \\ \hline
\end{tabular}
\end{center}

The quantity $Ai(x)^{2}$ satisf\thinspace ies the dif\thinspace ferential
equation 
\begin{equation*}
\frac{d\,^{3}Ai(x)^{2}}{d\,x^{\,3}}-4\,x\frac{d\,Ai(x)^{2}}{d\,x}%
-2\,Ai(x)^{2}=0,
\end{equation*}%
from which it follows that 
\begin{equation*}
\frac{d\,\,^{j+3}Ai(x)^{2}}{d\,x^{\,j+3}}-4\,x\frac{d^{\,\,j+1}\,Ai(x)^{2}}{%
d\,x^{\,j+1}}-(4j+2)\,\frac{d^{\,\,j}\,Ai(x)^{2}}{d\,x\,^{j}}=0.
\end{equation*}%
From this relation we get the recurrence equations 
\begin{eqnarray}
\mathfrak{P}_{\,j+3}(x)-4\,x\,\mathfrak{P}_{j+1}(x)-(4\,j+2)\,\mathfrak{P}%
_{j}(x) &=&0,  \notag \\
\mathfrak{Q}_{\,j+3}(x)-4\,x\,\mathfrak{Q}_{j+1}(x)-(4\,j+2)\,\mathfrak{Q}%
_{j}(x) &=&0,  \notag \\
\mathfrak{R\,}_{j+3}(x)-4\,x\,\mathfrak{R\,}_{j+1}(x)-(4\,j+2)\,\mathfrak{R\,%
}_{j}(x) &=&0.  \label{eq24}
\end{eqnarray}%
The generating functions $\Xi (t,x),\,\Lambda (t,x),$and $\varrho (t,x)$%
\begin{eqnarray*}
\Xi (t,x) &=&\sum_{j=0}^{\infty }\frac{t^{\,\,j}}{j!}\,\mathfrak{P}_{j}(x)\,,
\\
\Lambda (t,x) &=&\sum_{j=0}^{\infty }\frac{t^{\,\,j}}{j!}\,\mathfrak{Q}%
_{\,j}(x)\,, \\
\varrho (t,x) &=&\sum_{j=0}^{\infty }\frac{t^{\,\,j}}{j!}\,\mathfrak{R}%
\,_{j}(x),
\end{eqnarray*}%
can be obtained in closed form by computing their f\thinspace irst three
derivatives and using the recurrence relations (24) to obtain the
dif\thinspace ferential equations def\thinspace ining them. \ For example in
the case of $\ \Xi (t,x)$ we get 
\begin{equation*}
\frac{d^{\,3}\,\Xi (t,x)}{d\,t^{3}}-4(x+t)\frac{d^{\,}\,\Xi (t,x)}{d\,t}%
-2\,\Xi (t,x)=0,
\end{equation*}%
the solution of which is 
\begin{equation*}
\Xi (t,x)=g_{1}(x)\,Ai(t+x)^{2}+g_{2}\,(x)Bi(t+x)^{2}+g_{3}\,(x)\
Ai(t+x)\,Bi(t+x).
\end{equation*}%
The boundary conditions in this case are%
\begin{eqnarray*}
\lbrack \Xi (t,x)]_{t=0} &=&\mathfrak{P}_{0}(x)=1, \\
\lbrack d\,\Xi (t,x)/dt]_{t=0} &=&\mathfrak{P}_{1}(x)=0, \\
\lbrack d\,^{2}\Xi (t,x)/dt^{2}]_{t=0} &=&\mathfrak{P}_{2}(x)=2x,
\end{eqnarray*}%
with the resulting closed form expression for $\Xi (t,x)$ being 
\begin{equation*}
\Xi (t,x)=\pi ^{2}[Bi^{\,\prime }(x)Ai(x+t)-Ai^{\,\prime
}(x)Bi(x+t)\,]\,^{2}=\xi (t,x)^{2}.
\end{equation*}%
In a similar way, the closed form expressions for the generating functions $%
\Lambda (t,x)$ and $\varrho (t,x)$ are 
\begin{equation*}
\Lambda (t,x)=\pi ^{2}[Bi(x)Ai(x+t)\,-Ai(x)Bi(x+t)\,\,]^{2}=\lambda
(t,x)^{2},
\end{equation*}%
and%
\begin{eqnarray*}
\varrho (t,x) &=&-2\pi ^{2}[Bi(x)Bi^{\prime
}(x)\,Ai(x+t)^{2}+Ai(x)Ai^{\prime }(x)\,Bi(x+t)^{2}] \\
&&+2\pi ^{2}\{Ai(x)Bi^{\prime }(x)+Bi(x)Ai^{\prime }(x)\}\,Ai(x+t)\,Bi(x+t),
\end{eqnarray*}%
which can we written as%
\begin{eqnarray*}
\varrho (t,x) &=&2\pi ^{2}[Bi^{\,\prime }(z)\,Ai(z+t)-Ai^{\,\prime
}(z)\,Bi(z+t)] \\
&&\cdot \lbrack Ai(z)\,Bi(z+t)-Bi(z)\,Ai(z+t)] \\
&=&2\,\xi (t,x)\,\lambda (t,x).
\end{eqnarray*}%
As expected we see that $\Xi (t,x)$, $\Lambda (t,x)$ and $\varrho (t,x)$ are
related to the generating functions $\xi (t,x)$ and $\lambda (t,x)$
encountered above. \ We also note that the polynomials $\mathfrak{P}_{j}(x),%
\mathfrak{Q}_{\,j}(x)$ and $\mathfrak{R\,}_{j}(x)$ are related to the
polynomials $\mathcal{P}_{\,k}(z)$ and $Q\,_{k}(z)$ by simple bilinear
relations.

Expanding the generating functions $\Xi (-a,x),\,\Lambda (-a,x)$, and $%
\varrho (-a,x)$ in power series in $x$ we get an expression which contains
the incomplete Mellin transforms of $Ai(x)^{2},$ $Ai^{\,\prime }(x)^{2}$ and 
$Ai(x)\,Ai^{\,\prime }(x)$ 
\begin{eqnarray}
J_{n}(a) &=&\sum_{i=0}^{\infty }\tfrac{\Xi ^{(i)}(-a,0)}{i!}\int_{a}^{\infty
}x^{\,i-n}\,Ai(x)^{2}dx  \notag \\
&&+\sum_{i=0}^{\infty }\tfrac{\Lambda ^{(i)}(-a,0)}{i!}\int_{a}^{\infty
}x^{\,i-n}Ai^{\,\prime }(x)^{2}dx  \notag \\
&&+\,\sum_{i=0}^{\infty }\tfrac{\varrho \mathfrak{\,}^{(i)}(-a,0)}{i!}%
\int_{a}^{\infty }x^{\,i-n}Ai(x)\,Ai^{\,\prime }(x)dx.  \label{eq25}
\end{eqnarray}%
The f\thinspace irst few expansion coef\thinspace f\thinspace icients of
those appearing above are given by 
\begin{eqnarray*}
\Xi (-a,0) &=&\pi ^{2}Ai^{\,\prime }(0)^{2}\mathfrak{J}_{+}(a)^{2}, \\
\Xi ^{(1)}(-a,0) &=&2\pi ^{2}Ai^{\,\prime }(0)^{2}\mathfrak{J}_{+}(a)%
\mathfrak{J}_{+}^{\,\prime }(a), \\
\Xi ^{(2)}(-a,0) &=&2\pi ^{2}Ai^{\,\prime }(0)^{2}[\mathfrak{J}%
_{+}^{\,\prime }(a)^{2}-a\mathfrak{J}_{+}(a)^{2}\,] \\
&&-2\pi ^{2}Ai(0)Ai^{\,\prime }(0)\mathfrak{J}_{+}(a)\mathfrak{J}_{-}(a), \\
\Xi ^{(3)}(-a,0) &=&2\pi ^{2}Ai^{\,\prime }(0)^{2}[3\mathfrak{J}%
_{+}(a)^{2}-4a\,\mathfrak{J}_{+}(a)\mathfrak{J}_{+}^{\,\prime }(a)] \\
&&-6\pi ^{2}Ai(0)Ai^{\,\prime }(0)[\mathfrak{J}_{+}(a)\mathfrak{J}%
_{-}^{\,\prime }(a)+\mathfrak{J}_{+}^{\,\prime }(a)\mathfrak{J}_{-}(a)],
\end{eqnarray*}%
\begin{eqnarray*}
\Lambda (-a,0) &=&\pi ^{2}Ai(0)^{2}\mathfrak{J}_{-}(a)^{2}, \\
\Lambda ^{(1)}(-a,0) &=&2\pi ^{2}Ai(0)^{2}\mathfrak{J}_{-}(a)\mathfrak{J}%
_{-}^{\,\prime }(a)-2\pi ^{2}Ai(0)Ai^{\,\prime }(0)\mathfrak{J}_{+}(a)%
\mathfrak{J}_{-}(a), \\
\Lambda ^{(2)}(-a,0) &=&2\pi ^{2}Ai^{\,\prime }(0)^{2}\mathfrak{J}%
_{+}(a)^{2}+2\pi ^{2}Ai(0)^{2}[\mathfrak{J}_{-}^{\,\prime }(a)^{2}-a\,%
\mathfrak{J}_{-}(a)^{2}] \\
&&-4\pi ^{2}Ai(0)Ai^{\,\prime }(0)[\mathfrak{J}_{+}(a)\mathfrak{J}%
_{-}^{\,\prime }(a)+\mathfrak{J}_{+}^{\,\prime }(a)\mathfrak{J}_{-}(a)], \\
\Lambda ^{(3)}(-a,0) &=&12\pi ^{2}Ai^{\,\prime }(0)^{2}\mathfrak{J}_{+}(a)%
\mathfrak{J}_{+}^{\,\prime }(a) \\
&&+4\pi ^{2}Ai(0)^{2}[\mathfrak{J}_{-}(a)^{2}-2a\,\mathfrak{J}_{-}(a)%
\mathfrak{J}_{-}^{\,\prime }(a)] \\
&&+12\pi ^{2}Ai(0)Ai^{\,\prime }(0)\left[ a\,\mathfrak{J}_{+}(a)\mathfrak{J}%
_{-}(a)-\mathfrak{J}_{+}^{\,\prime }(a)\mathfrak{J}_{-}^{\,\prime }(a)\right]
,
\end{eqnarray*}

\begin{eqnarray*}
\varrho \,(-a,0) &=&2\pi ^{2}Ai(0)Ai^{\,\prime }(0)\mathfrak{J}_{+}(a)%
\mathfrak{J}_{-}(a), \\
\varrho ^{(1)}(-a,0) &=&-2\pi ^{2}Ai^{\,\prime }(0)^{2}\mathfrak{J}%
_{+}(a)^{2} \\
&&+2\pi ^{2}Ai(0)Ai^{\,\prime }(0)[\mathfrak{J}_{+}(a)\mathfrak{J}%
_{-}^{\,\prime }(a)-\mathfrak{J}_{+}^{\,\prime }(a)\mathfrak{J}_{-}(a)], \\
\varrho ^{(2)}(-a,0) &=&-8\pi ^{2}Ai^{\,\prime }(0)^{2}\mathfrak{J}_{+}(a)%
\mathfrak{J}_{+}^{\,\prime }(a)-2\pi ^{2}Ai(0)^{2}\mathfrak{J}_{-}(a)^{2} \\
&&+4\pi ^{2}Ai(0)Ai^{\,\prime }(0)[\mathfrak{J}_{+}^{\,\prime }(a)\mathfrak{J%
}_{-}^{\,\prime }(a)-a\,\mathfrak{J}_{+}(a)\,\mathfrak{J}_{-}(a)], \\
\varrho ^{(3)}(-a,0) &=&-12\pi ^{2}Ai^{\,}(0)^{2}\mathfrak{J}_{-}(a)%
\mathfrak{J}_{-}^{\,\prime }(a)+12\pi ^{2}Ai^{\,\prime }(0)^{2}[a\,\mathfrak{%
J}_{+}(a)^{2}-\mathfrak{J}_{+}^{\,\prime }(a)^{2}] \\
&&+8\pi ^{2}Ai(0)Ai^{\,\prime }(0)[2\mathfrak{J}_{+}(a)\mathfrak{J}_{-}(a) \\
&&-a\{\mathfrak{J}_{+}(a)\mathfrak{J}_{-}^{\,\prime }(a)+\mathfrak{J}%
_{+}^{\,\prime }(a)\mathfrak{J}_{-}(a)\}],
\end{eqnarray*}%
where the $\mathfrak{J}$ terms have been def\thinspace ined above.

\subsection{The incomplete Mellin transforms of $Ai(x)^{2},$ $Ai^{\,\prime
}(x)^{2},$ $Ai(x)\,Ai^{\,\prime }(x)$}

The incomplete Mellin transforms i.e. 
\begin{equation*}
\mathcal{I}_{n}(a)=\int_{a}^{\infty }x^{\,n}Ai(x)\,Ai^{\,\prime }(x)\,dx,
\end{equation*}%
and%
\begin{equation*}
i_{n}(a)=\int_{a}^{\infty }x^{\,n}Ai(x)^{2}\,dx,\;\;i_{n}^{\,\prime
}(a)=\int_{a}^{\infty }x^{\,n}Ai^{\,\prime }(x)^{2}\,dx,
\end{equation*}%
are analogues of the Mellin transforms $\ I_{n}(a),$ $I_{n}^{\,\prime }(a)$
which have appeared above. \ These transforms can also be written in closed
forms. \ Recurrence equations for these integrals have been given by Vallee 
\textit{et al }\cite{vallee} as 
\begin{equation*}
2(2n-1)\,\mathcal{I}_{n}(a)-n(n-1)(n-2)\,\mathcal{I}_{n-3}(a)=
\end{equation*}%
\begin{equation}
-(n-1)a^{n}Ai(a)^{2}-na^{n-1}Ai^{\,\prime
}(a)^{2}+n(n-1)a^{n-2}Ai(a)Ai^{\,\prime }(a),  \label{eq26}
\end{equation}%
with%
\begin{eqnarray*}
\mathcal{I}_{0}(a) &=&-\frac{1}{2}Ai(a)^{2}, \\
\mathcal{I}_{1}(a) &=&-\frac{1}{2}Ai^{\,\prime }(a)^{2}, \\
\mathcal{I}_{2}(a) &=&-\frac{a^{2}}{6}Ai(a)^{2}-\frac{a}{3}Ai^{\,\prime
}(a)^{2}+\frac{1}{3}Ai(a)Ai^{\,\prime }(a).
\end{eqnarray*}%
Solution of the dif\thinspace ference equation (26) yields\ (cf. Appendix
V)\ 
\begin{equation*}
\tfrac{12^{k+1}\Gamma (k+5/6)}{(3k)!}\mathcal{I}_{3k}(a)=Ai(a)^{2}%
\,B_{k,0}^{(0)}+Ai^{\,\prime }(a)^{2}\,B_{k,0}^{(1)}+Ai(a)\,Ai^{\,\prime
}(a)\,B_{k,0}^{(2)},
\end{equation*}%
\begin{equation*}
\tfrac{12^{k+1}\Gamma (k+7/6)}{(3k+1)!}\mathcal{I}_{3k+1}(a)=Ai(a)^{2}%
\,B_{k,1}^{(0)}+Ai^{\,\prime }(a)^{2}\,B_{k,1}^{(1)}+Ai(a)\,Ai^{\,\prime
}(a)\,B_{k,1}^{(2)},
\end{equation*}%
\begin{equation}
\tfrac{12^{k+1}\Gamma (k+3/2)}{(3k+2)!}\mathcal{I}_{3k+2}(a)=Ai(a)^{2}%
\,B_{k,2}^{(0)}+Ai^{\,\prime }(a)^{2}\,B_{k,2}^{(1)}+Ai(a)\,Ai^{\,\prime
}(a)\,B_{k,2}^{(2)},  \label{eq27}
\end{equation}%
where the $B$ polynomials are given by%
\begin{eqnarray*}
B_{k,0}^{(0)} &=&-\sum_{l=0}^{k}\tfrac{(3l-1)\,\Gamma
(l-1/6)\,(12a^{3})^{\,l}\,}{(3l)!}, \\
B_{k,0}^{(1)} &=&-\frac{1}{a}\sum_{l=1}^{k}\tfrac{\Gamma
(l-1/6)\,(12a^{3})^{\,l}\,}{(3l-1)!}, \\
B_{k,0}^{(2)} &=&\frac{1}{a^{2}}\sum_{l=1}^{k}\tfrac{\Gamma
(l-1/6)\,(12a^{3})^{\,l}\,}{(3l-2)!},
\end{eqnarray*}%
\begin{eqnarray*}
B_{k,1}^{(0)} &=&-a\sum_{l=1}^{k}\tfrac{(3\,l\,)\,\,\Gamma
(l+1/6)\,(12a^{3})^{\,l}}{(3l+1)!}, \\
B_{k,1}^{(1)} &=&-\sum_{l=0}^{k}\tfrac{\Gamma (l+1/6)\,(12a^{3})^{\,l}}{(3l)!%
}, \\
B_{k,1}^{(2)} &=&\frac{1}{a}\sum_{l=1}^{k}\tfrac{\,\Gamma
(l+1/6)\,(12a^{3})^{\,l}}{(3l-1)!},
\end{eqnarray*}%
\begin{eqnarray*}
B_{k,2}^{(0)} &=&-a^{2}\sum_{l=0}^{k}\tfrac{\,(3l+1)\,\Gamma
(l+1/2)\,(12a^{3})^{\,l}}{(3l+2)!}, \\
B_{k,2}^{(1)} &=&-a\sum_{l=0}^{k}\tfrac{\,\Gamma (l+1/2)\,(12a^{3})^{\,l}}{%
(3l+1)!}. \\
B_{k,2}^{(2)} &=&\sum_{l=0}^{k}\tfrac{\,\Gamma (l+1/2)\,(12a^{3})^{\,l}}{%
(3l)!}
\end{eqnarray*}%
The remaining integrals $\,i_{n}(a)$ and $i_{n}^{\,\prime }(a)$ are then
given in terms of the $\mathcal{I}_{n}(a)$ integrals by means of the Vallee
relations \cite{vallee} as%
\begin{eqnarray}
(2n+1)\,i_{n}(a) &=&-n(n-1)\mathcal{I}_{n-2}(a)  \notag \\
&&-a^{n+1}Ai(a)^{2}+a^{n}Ai^{\,\prime }(a)^{2}-na^{n-1}Ai(a)Ai^{\,\prime
}(a),  \label{eq28}
\end{eqnarray}%
\begin{eqnarray}
(2n+3)\,i_{n}^{\,\prime }(a) &=&-n(n+2)\mathcal{I}_{n-1}(a)+a^{%
\,n+2}Ai(a)^{2}  \notag \\
&&-a^{\,n+1}Ai^{\,\prime }(a)^{2}-(n+2)a^{\,n}Ai(a)Ai^{\,\prime }(a).
\label{eq29}
\end{eqnarray}%
Integration by parts \ of $i_{n}(a)$ and $i_{n}^{\,\prime }(a)$ gives the
simple forms ($n\neq -1$)%
\begin{eqnarray*}
\,i_{n}(a) &=&-\frac{1}{(n+1)}[2\,\mathcal{I}_{n+1}(a)+a^{\,n+1}Ai(a)^{2}],
\\
\,i_{n}^{\,\prime }(a) &=&-\frac{1}{(n+1)}[2\,\mathcal{I}_{n+2}(a)+a^{%
\,n+1}Ai^{\,\prime }(a)^{2}].
\end{eqnarray*}%
Elimination of $\mathcal{I}$ f\thinspace rom these two equations then gives 
\begin{equation*}
i_{n}^{\,\prime }(a)=\frac{1}{(n+1)}[(n+2)\,i_{n+1}(a)+a^{%
\,n+2}Ai(a)^{2}-a^{\,n+2}Ai^{\,\prime }(a)^{2}].
\end{equation*}

\subsection{The integrals $i_{-n}(a),$ $i_{-n}^{\,\,\prime }(a),$ and $%
\mathcal{I}_{-n}(a)$ ( $n>1$)}

The integrals $i_{-n}(a),$ $i_{-n}^{\,\,\prime }(a),$ and $\mathcal{I}%
_{-n}(a)$ are similarly interrelated and integration by parts of the
f\thinspace irst two of these integrals gives 
\begin{eqnarray*}
i_{-n}(a) &=&\tfrac{1}{(n-1)}[2\,\mathcal{I}_{-n+1}(a)+a^{-n+1}Ai(a)^{2}], \\
i_{-n}^{\,\,\prime }(a) &=&\tfrac{1}{(n-1)}[2\,\mathcal{I}%
_{-n+2}(a)+a^{-n+1}Ai^{\,\prime }(a)^{2}].
\end{eqnarray*}%
The recurrence relation for $\mathcal{I}_{\,-n\,}(a)$ is given by%
\begin{eqnarray*}
\mathcal{I}_{-n}(a) &=&\{n(n+1)(n+2)\mathcal{I}_{\,-n\,-3}(a) \\
&&-(n+1)Ai(a)^{2}/a^{n}-nAi^{\,\prime }(a)^{2}/a^{n+1} \\
&&-n(n+1)Ai(a)Ai^{\,\prime }(a)/a^{n+2}\}/(4n+2).
\end{eqnarray*}%
The solutions of this dif\thinspace ference equation are given by (30)%
\begin{eqnarray}
\tfrac{\Gamma (3k)}{12^{k-1}\Gamma (k+1/6)}\mathcal{I}_{-3k}(a) &=&\tfrac{12%
}{\Gamma (1/6)}\mathcal{I}_{\,-3}(a)+Ai(a)^{2}\,\beta _{k,0}^{(0)}  \notag \\
&&+Ai^{\,\prime }(a)^{2}\beta _{k,0}^{(1)}+Ai(a)\,Ai^{\,\prime }(a)\beta
_{k,0}^{(2)},  \notag \\
\tfrac{\Gamma (3k+1)}{12^{k-1}\Gamma (k+1/2)}\mathcal{I}_{-3k-1}(a) &=&%
\tfrac{12}{\sqrt{\pi }}\mathcal{I}_{\,-1}(a)+Ai(a)^{2}\,\beta _{k,1}^{(0)} 
\notag \\
&&+Ai^{\,\prime }(a)^{2}\beta _{k,1}^{(1)}+Ai(a)\,Ai^{\,\prime }(a)\beta
_{k,1}^{(2)},  \notag \\
\tfrac{\Gamma (3k+2)}{12^{k-1}\Gamma (k+5/6)}\mathcal{I}_{-3k-2}(a) &=&%
\tfrac{12}{\Gamma (5/6)}\mathcal{I}_{\,-2\,}(a)+Ai(a)^{2}\,\beta _{k,2}^{(0)}
\notag \\
&&+Ai^{\,\prime }(a)^{2}\beta _{k,2}^{(1)}+Ai(a)\,Ai^{\,\prime }(a)\beta
_{k,2}^{(2)},  \label{eq30}
\end{eqnarray}%
respectively where the $\beta $ polynomials are%
\begin{eqnarray*}
\beta _{k,0}^{(0)} &=&\sum_{l=0}^{k-1}\tfrac{(3l+1)\Gamma (3l)}{\Gamma
(l+7/6)(12a^{3})^{l}}, \\
\beta _{k,0}^{(1)} &=&\frac{1}{a}\sum_{l=0}^{k-1}\tfrac{\Gamma (3l)}{\Gamma
(l+7/6)(12a^{3})^{l}}, \\
\beta _{k,0}^{(2)} &=&\frac{1}{a^{2}}\sum_{l=0}^{k-1}\tfrac{\Gamma (3l+2)}{%
\Gamma (l+7/6)(12a^{3})^{l}},
\end{eqnarray*}%
\begin{eqnarray*}
\beta _{k,1}^{(0)} &=&\frac{1}{a}\sum_{l=0}^{k-1}\tfrac{(3l+2)\Gamma (3l+1)}{%
\Gamma (l+3/2)(12a^{3})^{l}}, \\
\beta _{k,1}^{(1)} &=&\frac{1}{a^{2}}\sum_{l=0}^{k-1}\tfrac{\Gamma (3l+1)}{%
\Gamma (l+3/2)(12a^{3})^{l}}, \\
\beta _{k,1}^{(2)} &=&\frac{1}{a^{3}}\sum_{l=0}^{k-1}\tfrac{\Gamma (3l+3)}{%
\Gamma (l+3/2)(12a^{3})^{l}},
\end{eqnarray*}

\begin{eqnarray*}
\beta _{k,2}^{(0)} &=&\frac{1}{a^{2}}\sum_{l=0}^{k-1}\tfrac{(3l+3)\Gamma
(3l+2)}{\Gamma (l+11/6)(12a^{3})^{l}}, \\
\beta _{k,2}^{(1)} &=&\frac{1}{a^{3}}\sum_{l=0}^{k-1}\tfrac{\Gamma (3l+2)}{%
\Gamma (l+11/6)(12a^{3})^{l}}, \\
\beta _{k,2}^{(2)} &=&\frac{1}{a^{4}}\sum_{l=0}^{k-1}\tfrac{\Gamma (3l+4)}{%
\Gamma (l+11/6)(12a^{3})^{l}}.
\end{eqnarray*}%
Written in these terms the expression for $J_{n}(a)$ is 
\begin{eqnarray}
J_{n}(a) &=&\sum_{k=0}^{\infty }\{\tfrac{\Xi ^{(k+n)}}{(k+n)!}\,\,i_{k}(a)+%
\tfrac{\Lambda ^{(k+n)}}{(k+n)!}\,i_{k}^{\,\,\prime }(a)+\,\tfrac{\mathfrak{%
\,\varrho }^{(k+n)}}{(k+n)!}\,\mathcal{I}_{k}(a)\}  \notag \\
&&+\sum_{k=1}^{n}\{\tfrac{\Xi ^{(n-k)}}{(n-k)!}\,i_{-k}(a)+\tfrac{\Lambda
^{(n-k)}}{(n-k)!}\,i_{-k}^{\,\,\prime }(a)+\,\tfrac{\mathfrak{\varrho }%
^{(n-k)}}{(n-k)!}\,\mathcal{I}_{-k}(a)\}.  \label{eq31}
\end{eqnarray}%
A general expression for $J_{n}(a)$ in its most reduced form for all $n$ is
quite complicated and is given in Appendix VI.

Integrals in (31) which contain negative powers of $x$ such as 
\begin{eqnarray*}
i_{-1} &(a)=&\int_{a}^{\infty }\frac{Ai(x)^{2}}{x}dx, \\
i_{-1}^{\,\prime } &(a)=&\int_{a}^{\infty }\frac{Ai\,^{\prime }(x)^{2}}{x}dx,
\\
\mathcal{I}_{-1}(a) &=&\int_{a}^{\infty }\frac{Ai(x)\,Ai^{\,\,\prime }(x)}{x}%
dx,
\end{eqnarray*}%
are irreducible and require special treatment (integrals with larger
negative powers of $x$ being expressible in terms of those integrals). \ We
have 
\begin{eqnarray*}
i_{-1}(a) &=&-Ai(0)^{2}[\,\ln (a)+\tfrac{1}{9}a^{3}\,_{3}F_{4}\tbinom{%
1,1,7/6\;;\;4a^{3}/9}{4/3,5/3,2,2}] \\
&&-2\,a\,Ai(0)\,Ai^{\,\prime }(0)\,\,_{2}F_{3}\tbinom{1/3,1/2\;;\;4a^{3}/9}{%
2/3,4/3,4/3} \\
&&-\tfrac{1}{2}a^{\,2}\,Ai^{\,\prime }(0)^{2}\,_{2}F_{3}\tbinom{%
2/3,5/6\;;\;4a^{3}/9}{4/3,5/3,5/3},
\end{eqnarray*}%
\begin{eqnarray*}
i_{-1}^{\,\prime }(a) &=&-Ai^{\,\prime }(0)^{2}[\,\ln (a)+\tfrac{2}{9}%
a^{3}\,_{3}F_{4}\tbinom{5/6,1,1\;;\;4a^{3}/9}{2/3,4/3,2,2}] \\
&&-\tfrac{1}{2}a^{2}Ai(0)\,Ai^{\,\prime }(0)\,\,_{2}F_{3}\tbinom{%
1/2,2/3\;;\;4a^{3}/9}{1/3,5/3,5/3} \\
&&-\tfrac{1}{16}a^{4}Ai(0)^{2}\,\,_{2}F_{3}\tbinom{7/3,4/3\;;\;4a^{3}/9}{%
5/3,7/3,7/3},
\end{eqnarray*}%
\begin{eqnarray*}
\mathcal{I}_{-1}(a) &=&-\tfrac{1}{4}a^{2}\,Ai(0)^{2}{}_{2}\,F_{3}\,\tbinom{%
2/3,7/6\;;\;4a^{3}/9}{4/3,5/3,5/3} \\
&&-Ai(0)\,Ai^{\,\prime }(0)[\,\ln (a)+\tfrac{1}{3}a^{3}\,_{3}F_{4}\tbinom{%
1,1,3/2\;;\;4a^{3}/9}{4/3,5/3,5/3,2}] \\
&&-a\,\,Ai^{\,\prime }(0)^{2}\,_{2}F_{3}\tbinom{1/3,5/6\;;\;4a^{3}/9}{%
2/3,4/3,4/3},
\end{eqnarray*}%
and were obtained by integrating the Taylor series representations for the
Airy products and then summing those series. \ Mathematica gives the same
expressions for these integrals. \ Some additional integrals which contain
negative indices $n$ are given by%
\begin{eqnarray*}
i_{-2}(a) &=&\tfrac{1}{a}Ai(a)^{2}+2\,\mathcal{I}_{-1}(a), \\
i_{-3}(a) &=&Ai(a)^{2}(\tfrac{1}{2a^{2}}-a)+Ai^{\,\,\prime }(a)^{2}+\tfrac{1%
}{a}Ai(x)Ai^{\,\,\prime }(a)+i_{-1}^{\,\,\prime }(a), \\
i_{-2}^{\,\,\prime }(a) &=&-Ai(a)^{2}+\tfrac{1}{a}Ai^{\,\,\prime }(a)^{2}, \\
i_{-3}^{\,\,\prime }(a) &=&\tfrac{1}{2a^{2}}Ai^{\,\,\prime }(a)^{2}+\,%
\mathcal{I}_{-1}(a), \\
\mathcal{I}_{-2}(a) &=&-a\,Ai(a)^{2}+Ai^{\,\,\prime }(a)^{2}+\tfrac{1}{a}%
Ai(a)Ai^{\,\,\prime }(a)+i_{-1}^{\,\,\prime }(a), \\
\mathcal{I}_{-3}(a) &=&-\tfrac{1}{2}Ai^{\,}(a)^{2}+\tfrac{1}{2a}%
Ai^{\,\,\prime }(a)^{2}+\tfrac{1}{2a^{2}}Ai(x)Ai^{\,\,\prime }(a)+\tfrac{1}{2%
}i_{-1}(a).
\end{eqnarray*}

We note in passing that in the case of large $a,$ the former integrals have
the asymptotic forms 
\begin{eqnarray*}
i_{-1}(a) &=&\int_{a}^{\infty }\frac{Ai(x)^{2}}{x}dx\thicksim \frac{e^{-%
\frac{4}{3}a^{3/2}}}{8\pi a^{2}}, \\
i_{-1}^{\,\prime }(a) &=&\int_{a}^{\infty }\frac{Ai^{\,\prime }(x)^{2}}{x}%
dx\thicksim \frac{e^{-\frac{4}{3}a^{3/2}}}{8\pi a}, \\
\mathcal{I}_{-1}(a) &=&\int_{a}^{\infty }\frac{Ai(x)\,Ai^{\,\prime }(x)}{x}%
dx\thicksim -\frac{e^{-\frac{4}{3}a^{3/2}}}{8\pi a^{3/2}}.
\end{eqnarray*}%
where the asymptotic expressions for the Airy functions have been used
above. \ \bigskip

Integrations based upon the expressions given in Appendix VI, for $%
J_{1}(\left| a_{1}^{\prime }\right| ),$ $J_{2}(\left| a_{1}^{\prime }\right|
),$ and $J_{3}(\left| a_{1}^{\prime }\right| )$ yield values of $0.04826441,$
$0.03654795,$and $0.02879281$ respectively. \ Values accurate to eleven
decimal places when powers of $x$ up to $x^{8}$ are included.\ 

\subsection{Conclusions}

Although the integrals $\mathfrak{I}^{\,\left( 1\right) },$ and$\ \mathfrak{I%
}^{\,\left( 2\right) }$ have been given in terms of analytic functions, it
is unfortunate that these expressions are extremely complicated in form and
thus of limited usefulness. This is disappointing since the Thomas-Fermi
formulation for atomic theory is in essence an analytic i.e. non-numeric
approach such as Hartree-Fock methods used in the study of atomic species
and analytic evaluation of its integrals would be in keeping with the spirit
of the method. \ \pagebreak

\begin{center}
\appendix Appendix I
\end{center}

Using the Weirstrass representation for the quantity $Ai^{\,\prime
}(z)/Ai^{\,\prime }(0)$ we have 
\begin{eqnarray*}
\tfrac{d\,\ln (Ai^{\,\prime }(z)/Ai^{\,\prime }(0))}{d\text{\/}\,z} &=&%
\tfrac{z\,Ai(z)}{Ai^{\,\prime }(z)}, \\
&=&\sum_{n=1}^{\infty }[\tfrac{1}{(z+\left| a_{n}^{\,\prime }\right| )}-%
\tfrac{1}{\left| a_{n}^{\,\prime }\right| }].
\end{eqnarray*}%
Then for $k>1$ we can obtain from the equation above the Airy zeta function $%
\mathcal{Z}_{k}$ which is def\thinspace ined as%
\begin{equation*}
\mathcal{Z}_{k}\,\mathcal{=}\sum_{n=1}^{\infty }\frac{1}{\left|
a_{n}^{\,\prime }\right| ^{\,k}},
\end{equation*}%
by means of the expression 
\begin{equation*}
\mathcal{Z}_{k}=\frac{(-1)^{k-1}}{(k-1)!}[\frac{d^{\,k-1}}{d\,z^{k-1}}\{%
\frac{zAi(z)}{Ai^{\,\prime }(z)}\}]_{z=0}.
\end{equation*}%
\ The incomplete Airy zeta function $\mathcal{Z}_{k}(N)$ is def\thinspace
ined as the f\thinspace inite sum%
\begin{equation*}
\mathcal{Z}_{k}(N)\,\mathcal{=}\sum_{n=1}^{N}\frac{1}{\left| a_{n}^{\,\prime
}\right| ^{\,k}}\,.
\end{equation*}%
The f\thinspace irst few of the $\mathcal{Z}_{k}$ sums are given in\ Table 5
( $\eta =Ai(0)/Ai^{\,\,\prime }(0)$).

\begin{center}
\begin{table}[H] \centering%
\caption{The Airy zeta function Z\label{six}}%
\end{table}%

\vspace{0.5in}

\renewcommand{\arraystretch}{2.0}%
\begin{tabular}{|c|c|}
\hline
$k$ & $\mathcal{Z}_{k}$ \\ \hline
$2$ & $\eta $ \\ \hline
$3$ & $1$ \\ \hline
$4$ & $\frac{1}{2}\eta ^{2}$ \\ \hline
$5$ & $-\frac{2}{3}\eta $ \\ \hline
$6$ & $\frac{1}{4}-\frac{1}{4}\eta ^{3}$ \\ \hline
$7$ & $\frac{7}{15}\eta ^{2}$ \\ \hline
$8$ & $-\frac{11}{36}\eta +\frac{1}{8}\eta ^{4}$ \\ \hline
\end{tabular}

\bigskip

\pagebreak

\appendix Appendix II
\end{center}

It is possible to obtain the leading terms in an asymptotic expression for
the integral $\mathbf{I}_{1}(a)$ as follows. \ Rewriting the integral as%
\begin{equation*}
\mathbf{I}_{1}(a)=\frac{1}{a}\int_{0}^{\infty }\frac{Ai(z)}{(1+z/a)}dz,
\end{equation*}%
expansion of the denominator for large $a$ yields 
\begin{eqnarray*}
\mathbf{I}_{1}(a) &\thicksim &\frac{1}{a}\sum_{k=0}^{\infty }\left[ \frac{-1%
}{a}\right] ^{k}\int_{0}^{\infty }z^{\,k}Ai(z)\text{\/\thinspace }dz, \\
\mathbf{I}_{1}(a) &\thicksim &\frac{1}{3a}\sum_{k=0}^{\infty }\left[ \frac{-1%
}{3^{1/3}a}\right] ^{k}\frac{k!}{\Gamma (k/3+1)}, \\
\mathbf{I}_{1}(a) &\thicksim &\frac{1}{3a}-\frac{1}{3^{1/3}\,\Gamma
(1/3)\,a^{2}}+\cdots
\end{eqnarray*}%
Repeated dif\thinspace ferentiation of the expression for $\mathbf{I}_{1}(a)$
then produces large $a$ expressions for $\mathbf{I}_{k}(a).$

\begin{center}
\bigskip

\appendix Appendix III
\end{center}

We note that the integrals $I_{n}^{^{\;}}(a)$ and $I_{n}^{^{\;\prime }}(a)$
have for $n>0$ the forms%
\begin{eqnarray*}
I_{n}(a) &=&\mathrm{c}_{n}(a)\,Ai(a)+\mathrm{d}_{n}(a)\,Ai^{\,\prime }(a)+%
\mathrm{e}_{n}\,I_{0}(a), \\
I_{n}^{\,\,\prime }(a) &=&\mathrm{c}_{n}^{\,\prime }(a)\,Ai(a)+\mathrm{d\,}%
_{n}^{\prime }(a)\,Ai^{\,\prime }(a)+\mathrm{e}_{n}^{\,\prime }\,I_{0}(a).
\end{eqnarray*}%
Using the recurrence relation (15) we f\thinspace ind that $\mathrm{c}%
_{n}(a),$ $\mathrm{d}_{n}(a),$ are polynomials which satisfy the relations%
\begin{eqnarray*}
\mathrm{c}_{n}(a) &=&(n-1)(n-2)\,\mathrm{c}_{n-3}(a)+(n-1)\,a^{n-2}, \\
\mathrm{d}_{n}(a) &=&(n-1)(n-2)\,\mathrm{d}_{n-3}(a)-a^{n-1},
\end{eqnarray*}%
and the constants $\mathrm{e}_{n}$ are given by%
\begin{equation*}
\mathrm{e}_{n}=(n-1)(n-2)\,\mathrm{e}_{n-3}.
\end{equation*}%
Initial values of these quantities are given in the f\thinspace irst three
rows of Table 6 along with their f\thinspace irst few values.

\begin{center}
\begin{table}[H] \centering%
\caption{The polynomials c,d and e\label{seven}}%
\end{table}%

\begin{tabular}{|c|c|c|c|}
\hline
$n$ & $\mathrm{c}_{n}(a)$ & $\mathrm{d}_{n}(a)$ & $\mathrm{e}_{n}$ \\ \hline
$1$ & $0$ & $-1$ & $0$ \\ \hline
$2$ & $1$ & $-a$ & $0$ \\ \hline
$3$ & $2a$ & $-a^{2}$ & $2$ \\ \hline
$4$ & $3a^{2}$ & $-(6+a^{3})$ & $0$ \\ \hline
$5$ & $12+4a^{3}$ & $-(12a+a^{4})$ & $0$ \\ \hline
$6$ & $40+5a^{4}$ & $-(20a^{2}+a^{5)}$ & $40$ \\ \hline
\end{tabular}
\end{center}

The solution for $\mathrm{e}_{n}$ follows immediately and we get%
\begin{equation*}
\mathrm{e}_{n}=[1+2\cos (2\pi n/3)]\frac{(n-1)!}{3^{n/3-2}\Gamma (n/3)},
\end{equation*}%
which has non-zero value only for the case $\mathrm{e}_{3k}(a)$ i.e.%
\begin{equation*}
\mathrm{e}_{3k}=\frac{(3k)!}{3^{k}\,k!},\hspace{0.25in}\mathrm{e}_{3k+1}=0,%
\hspace{0.25in}\mathrm{e}_{3k+2}=0.
\end{equation*}%
In addition we have 
\begin{eqnarray*}
\mathrm{c}_{3k}(a) &=&\left( \frac{a}{3}\right) 9^{k}\Gamma (k+1/3)\Gamma
(k+2/3)\sum_{l=0}^{k-1}\frac{(a^{3}/9)^{l}}{\Gamma (l+2/3)\,\Gamma (l+4/3)},
\\
\mathrm{c}_{3k+1}(a) &=&\left( \frac{a^{2}}{3}\right) 9^{k}k!\,\Gamma
(k+2/3)\sum_{l=0}^{k-1}\frac{(a^{3}/9)^{l}}{l!\,\Gamma (l+5/3)}, \\
\mathrm{c}_{3k+2}(a) &=&3\cdot 9^{k}k!\Gamma (k+4/3)\sum_{l=0}^{k}\frac{%
(a^{3}/9)^{l}}{l!\,\Gamma (l+1/3)},
\end{eqnarray*}%
\begin{eqnarray*}
\mathrm{d}\,_{3k}(a) &=&-\left( \frac{a^{2}}{9}\right) 9^{k}\Gamma
(k+1/3)\Gamma (k+2/3)\sum_{l=0}^{k-1}\frac{(a^{3}/9)^{l}}{\Gamma
(l+4/3)\Gamma (l+5/3)}, \\
\mathrm{d}\,_{3k+1}(a) &=&-9^{k}k!\Gamma (k+2/3)\sum_{l=0}^{k}\frac{%
(a^{3}/9)^{l}}{l!\,\Gamma (l+2/3)}, \\
\mathrm{d}\,_{3k+2}(a) &=&-a\,9^{k}k!\Gamma (k+4/3)\sum_{l=0}^{k}\frac{%
(a^{3}/9)^{l}}{l!\,\Gamma (l+4/3)}.
\end{eqnarray*}%
The integrals $I_{n}^{\prime }(a)$ in terms of these polynomials are 
\begin{equation*}
I_{n}^{\prime }(a)=-[a^{n}+n\,\mathrm{c}_{n-1}(a)]\,Ai(a)-n\,\mathrm{d}%
\,_{n-1}(a)\,Ai^{\,\,\prime }(a)-n\,\mathrm{e}_{n-1}\,I_{0}(a)\,.
\end{equation*}

\subsection{\protect\Large The Mellin transforms for the Airy products}

\bigskip

Mellin transforms of the functions $Ai(x)^{2},$ $Ai^{\,\prime }(x)^{2},$ $%
Ai(x)\,Ai^{\,\prime }(x)$ may be obtained from well known recurrence
relations. \ For example, the moment $\widehat{i}_{n}$ has the recurrence
relation 
\begin{equation*}
\widehat{i}_{j}(0)=\int_{0}^{\infty }x^{\,j}Ai(x)Ai^{\,\prime }(x)\,dx=%
\tfrac{j\,(j-1)(j-2)}{2(2j-1)}\int_{0}^{\infty }x^{\,j-3}Ai(x)Ai^{\,\prime
}(x)\,dx,
\end{equation*}%
or viewed as a dif\thinspace f\thinspace erence equation i.e.%
\begin{equation*}
\widehat{i}_{j}=\frac{j\,(j-1)(j-2)}{2(2j-1)}\,\widehat{\,i}_{j-3},
\end{equation*}%
with initial values 
\begin{eqnarray*}
\widehat{i}_{0}(0) &=&-\tfrac{1}{2}Ai(0)^{2}, \\
\widehat{i}_{1}(0) &=&-\tfrac{1}{2}Ai^{\prime }(0)^{2}, \\
\widehat{i}_{2}(0) &=&\tfrac{1}{3}Ai(0)\,Ai^{\prime }(0),
\end{eqnarray*}%
the $\widehat{i}_{j}$ integrals have values given by 
\begin{eqnarray*}
\widehat{i}_{3j}(0) &=&-\frac{(3j)!\,\Gamma (5/6)}{2\,(12)^{\,\,j}\,\,\Gamma
(j+5/6)}\,Ai(0)^{2}, \\
\widehat{i}_{3j+1}(0) &=&-\frac{(3j+1)!\,\pi }{6\,(12)^{\,\,j}\,\Gamma
(5/6)\Gamma (j+7/6)}\,Ai^{\,\prime }(0)^{2}, \\
\widehat{i}_{3j+2}(0) &=&\frac{(3j+2)!\,\sqrt{\pi }}{(12)^{\,\,j+1}\,\Gamma
(j+3/2)}\,Ai(0)Ai^{\,\prime }(0).
\end{eqnarray*}%
The moments of $Ai(x)^{2},$and $Ai^{\,\prime }(x)^{2}$ are then seen to be
multiples of the $\widehat{i}_{j}(0)$ integrals. \ We have%
\begin{eqnarray*}
\int_{0}^{\infty }x^{\,j}Ai(x)^{2}dx &=&-\frac{j\,(j\,-1)}{2\,j+1}\,\widehat{%
i}_{j-2}(0), \\
\int_{0}^{\infty }x^{\,j}Ai^{\,\prime }(x)^{2}dx &=&-\frac{j\,(j+2)}{2\,j+3}%
\,\widehat{i}_{j-1}(0).
\end{eqnarray*}%
More explicitly%
\begin{eqnarray*}
\int_{0}^{\infty }x^{3j+2}Ai(x)^{2}dx &=&\frac{(3j+2)!\,\Gamma (5/6)}{%
(12)^{\,j+1}\,\Gamma (j+11/6)}\,Ai(0)^{2}, \\
\int_{0}^{\infty }x^{3j+3}Ai(x)^{2}dx &=&\frac{\pi \,(3j+3)!}{%
3\,(12)^{\,j+1}\,\Gamma (5/6)\Gamma (j+13/6)}\,Ai^{\,\prime \,}(0)^{2}, \\
\int_{0}^{\infty }x^{3j+4}Ai(x)^{2}dx &=&-\frac{2\sqrt{\pi }\,(3j+4)!}{%
(12)^{\,j+2}\,\Gamma (j+5/2)}\,Ai(0)\,Ai^{\,\prime }(0),
\end{eqnarray*}%
and%
\begin{eqnarray*}
\int_{0}^{\infty }x^{3j+1}Ai^{\,\prime }(x)^{2}dx &=&\frac{%
\,(3j+3)\,(3j+1)!\,\Gamma (11/6)}{10\,(12)^{\,j}\,\Gamma (j+11/6)}%
\,Ai(0)^{2}, \\
\int_{0}^{\infty }x^{3j+2}Ai^{\,\prime }(x)^{2}dx &=&\frac{5\,\pi
\,(3j+4)\,(3j+2)!\,}{18\,(12)^{\,j+1}\,\Gamma (11/6)\,\Gamma (j+13/6)}%
\,Ai^{\,\prime }(0)^{2}, \\
\int_{0}^{\infty }x^{3j+3}Ai^{\,\prime }(x)^{2}dx &=&-\frac{2\sqrt{\pi }%
\,(3\,j+5)\,(3j+3)!}{\,(12)^{\,j+2}\,\Gamma (j+5/2)}\,Ai(0)Ai^{\,\prime }(0).
\end{eqnarray*}%
The integrals above can be written in more compact form using relations
obtained by\ Reid \cite{Reid} i.e.%
\begin{eqnarray*}
\int_{0}^{\infty }x^{\,\alpha -1}Ai(x)^{2}dx &=&\frac{2\,\Gamma (\alpha )}{%
\sqrt{\pi }12^{\,(\alpha /3+5/6)}\Gamma (\alpha /3+5/6)}, \\
\int_{0}^{\infty }x^{\,\alpha -1}Ai^{\,\prime }(x)^{2}dx &=&\frac{2(\alpha
+1)\,\Gamma (\alpha )}{\sqrt{\pi }12^{\,(\alpha /3+7/6)}\Gamma (\alpha
/3+7/6)}, \\
\int_{0}^{\infty }x^{\,\alpha -1}Ai(x)\,Ai^{\,\prime }(x)\,dx &=&-\,\frac{%
2(2\alpha +3)\,\Gamma (\alpha )}{\sqrt{\pi }12^{\,(\alpha /3+3/2)}\Gamma
(\alpha /3+3/2)}.
\end{eqnarray*}

\pagebreak \bigskip

\begin{center}
\appendix Appendix IV
\end{center}

A general expression for $\mathbf{I}_{n}(a)$ 
\begin{equation*}
\mathbf{I}_{n}(a)=-\frac{Ai(a)}{a^{n}}\sum_{i=0}^{\infty }\frac{%
a^{\,i}\lambda ^{(i)}}{i!}+\sum_{i=0}^{\infty
}b_{i}\,I_{i}(a)+\sum_{i=1}^{n+1}b_{-i\,}I_{-\,i}(a).
\end{equation*}%
can be written in terms of only the $I_{\,\pm \,i\,}(a)$ integrals as
follows. The f\thinspace irst sum above can be written in closed form as%
\begin{equation*}
\sum_{i=0}^{\infty }\frac{a^{\,i}\,\lambda ^{(i)}}{i!}=\lambda (-a,a)=\pi
Ai(0)\,\mathfrak{J}_{-}(-a),
\end{equation*}%
and the coef\thinspace f\thinspace icients $b_{i}$ are given by%
\begin{equation*}
b_{i}=\frac{1}{(n+1+i)!}\left\{ (n+1+i)\xi ^{(n+i)}-(1+i)\lambda
^{(n+1+i)}\right\} .
\end{equation*}%
The last two sums in $\mathbf{I}_{n}(a)$ can be written as 
\begin{eqnarray*}
\sum_{i=0}^{\infty }b_{i}\,I_{i}(a) &=&\sum_{i=0}^{\infty
}b_{3i}\,I_{3i}(a)+\sum_{i=0}^{\infty }b_{3i+1}\,I_{3i+1}(a) \\
&&+\sum_{i=0}^{\infty }b_{3i+2}\,I_{3i+2}(a), \\
\sum_{i=1}^{n+1}b_{-i\,}I_{-\,i}(a) &=&\sum_{i=1}^{\lfloor \frac{n+1}{3}%
\rfloor }b_{-3i\,}I_{-3\,i}(a)+\sum_{i=0}^{\lfloor \frac{n}{3}\rfloor
}b_{-3i-1\,}I_{-3\,i-1}(a) \\
&&+\sum_{i=0}^{\lfloor \frac{n-1}{3}\rfloor }b_{-3i\,-2}I_{-3\,i-2}(a),
\end{eqnarray*}%
respectively, where $\lfloor z\rfloor $ is the f\thinspace loor function. \
Gathering together the terms common to $I_{0}(a),\,Ai(a)$ and $Ai^{\,\prime
}(a)$ we have%
\begin{eqnarray*}
\sum_{i=0}^{\infty }b_{i}I_{i}(a) &=&I_{0}(a)\sum_{k=0}^{\infty }\Omega
_{k},_{0}(n|a)+Ai(a)\sum_{\mu =0}^{2}\sum_{k=0}^{\infty }\Omega _{k},_{\mu
}(n|a)\,\sigma _{k,\,\mu }(a) \\
&&-Ai^{\,\prime }(a)\sum_{\mu =0}^{2}\sum_{k=0}^{\infty }\Omega _{k},_{\mu
}(n|a)\,\sigma _{k,\,\mu }^{\prime }(a),
\end{eqnarray*}%
where 
\begin{equation*}
\Omega _{k},_{\mu }(n|a)=\tfrac{1}{3^{k+\mu }(3k+\mu +2)_{n-1}\Gamma (k+1+%
\frac{\mu }{3})}[\frac{\xi ^{(3k+n+\mu )}\,}{(3k+\mu +1)}-\,\frac{\lambda
^{(3k+n+\mu +1)}\,}{(3k+n+\mu +1)}],
\end{equation*}%
and $(z)_{n}$ is the Pochhammer function. The polynomials $\sigma _{k},_{\mu
}$ and $\sigma _{k}^{\,\prime },_{\mu }$ appearing above are given by%
\begin{eqnarray*}
\sigma _{k,\,0}(a) &=&\,a\sum_{l=0}^{k-1}\tfrac{(3a^{3})^{l}\,l!}{(3l+1)!},
\\
\sigma _{k,\,1}(a) &=&a^{2}\,\frac{2\pi }{3\sqrt{3}}\sum_{l=0}^{k-1}\tfrac{%
(a^{3}/9)^{l}\,}{l\,!\,\,\Gamma (l+5/3)}, \\
\sigma _{k,\,2}(a) &=&2\pi \sqrt{3}\sum_{l=0}^{k}\tfrac{(a^{3}/9)^{l}\,}{%
l\,!\,\,\,\Gamma (l+1/3)},
\end{eqnarray*}%
and%
\begin{eqnarray*}
\sigma _{k,\,0}^{\,\prime }(a) &=&a^{2}\sum_{l=0}^{k-1}\tfrac{%
(3a^{3})^{l}\,l!}{(3l+2)!}, \\
\sigma _{k,1}^{\,\prime }(a) &=&\frac{2\pi }{\sqrt{3}}\sum_{l=0}^{k}\tfrac{%
(a^{3}/9)^{l}\,}{l\,!\,\,\Gamma (l+2/3)}, \\
\sigma _{k,2}^{\,\prime }(a) &=&a\,\frac{2\pi }{\sqrt{3}}\sum_{l=0}^{k}%
\tfrac{(a^{3}/9)^{l}\,}{l\,!\,\,\Gamma (l+4/3)}.
\end{eqnarray*}%
The f\thinspace inite sum can be rewritten in a similar way as 
\begin{eqnarray*}
\sum_{i=1}^{n+1}b_{-i}\,I_{-i}(a) &=&I_{0}(a)\{\tfrac{2\pi \emph{\ }}{\sqrt{3%
}}\sum_{i=1}^{\lfloor \frac{n+1}{3}\rfloor }\omega _{i,0}(n|a)\} \\
&&+I_{-1}(a)\Gamma (2/3)\sum_{i=0}^{\lfloor \frac{n}{3}\rfloor }\omega
_{i,1}(n|a) \\
&&+I_{-2}(a)\tfrac{1}{3}\Gamma (1/3)\sum_{i=0}^{\lfloor \frac{n-1}{3}\rfloor
}\omega _{i,2}(n|a)\emph{\ \ \ } \\
&&+Ai(a)\sum_{\mu =\,0}^{2}\sum_{i=1}^{\lfloor \frac{n+1-\mu }{3}\rfloor
}\omega _{i,\mu }(n|a)\,\,\widehat{\sigma }_{i,\mu }(a) \\
&&+Ai^{\,\prime }(a)\sum_{\mu =\,0}^{2}\sum_{i=1}^{\lfloor \frac{n+1-\mu }{3}%
\rfloor }\omega _{i,\mu }(n|a)\,\,\widehat{\sigma }_{i,\mu }^{\,\prime }(a),%
\emph{\ }
\end{eqnarray*}%
where 
\begin{equation*}
\omega _{\,i,\,\mu }(n|a)=\tfrac{1}{3^{2i-1}(n-\mu -3\,i)\,!\,\Gamma (i+%
\frac{\mu +1}{3})\Gamma (i+\frac{\mu -1}{3})}[\tfrac{\xi ^{(n-\mu -3i)}\,}{%
(3i+\mu -1)}+\tfrac{\lambda ^{(n+1-\mu -3i)}\,}{(n+1-\mu -3i)}],
\end{equation*}%
(note the cases $\omega _{0,2}=\frac{\xi ^{(2)}}{\Gamma (2/3)},$ $\omega
_{1,2}=\frac{\lambda ^{(0)}}{2\Gamma (2/3)}$). \ The polynomials $\,\widehat{%
\sigma }_{k,\,\mu }$ ,$\,\widehat{\sigma }_{k,\,\mu }^{\,\prime }$ are 
\begin{eqnarray*}
\widehat{\sigma }_{k,0}(a) &=&\frac{2\pi }{a^{2}\sqrt{3}}\sum_{l=0}^{k-1}%
\frac{(3l+1)!}{l!\,\left( 3a^{3}\right) ^{l}}, \\
\widehat{\sigma }_{k,1}(a) &=&\frac{3}{a^{3}}\sum_{l=0}^{k-1}\frac{%
l!\,\Gamma (l+5/3)}{\left( a^{3}/9\right) ^{l}}, \\
\widehat{\sigma }_{k,2}(a) &=&\frac{3}{a^{4}}\sum_{l=0}^{k-1}\frac{%
(l+1)!\,\Gamma (l+4/3)}{\left( a^{3}/9\right) ^{l}},
\end{eqnarray*}%
and%
\begin{eqnarray*}
\widehat{\sigma }_{k,0}^{\,\prime }(a) &=&\frac{2\pi }{a\sqrt{3}}%
\sum_{l=0}^{k-1}\frac{(3l)!}{l!\,\left( 3a^{3}\right) ^{l}}, \\
\widehat{\sigma }_{k,1}^{\,\prime }(a) &=&\frac{1}{a^{2}}\sum_{l=0}^{k-1}%
\frac{l!\,\Gamma (l+2/3)}{\left( a^{3}/9\right) ^{l}}, \\
\widehat{\sigma }_{k,2}^{\,\prime }(a) &=&\frac{1}{a^{3}}\sum_{l=0}^{k-1}%
\frac{l!\,\Gamma (l+4/3)}{\left( a^{3}/9\right) ^{l}}.
\end{eqnarray*}%
\ Finally we get for $\mathbf{I}_{n}(a)$%
\begin{eqnarray*}
\mathbf{I}_{n}(a) &=&-\frac{\pi Ai(0)Ai(a)}{a^{n}}\mathfrak{J}_{-}(-a)+ \\
&&I_{0}(a)\{\sum_{k=0}^{\infty }\Omega _{k,0}(n|a)+\tfrac{2\pi }{\sqrt{3}}%
\sum_{k=1}^{\lfloor \frac{n+1}{3}\rfloor }\omega _{i,0}(n|a)\} \\
&&+I_{-1}(a)\Gamma (2/3)\sum_{k=1}^{\lfloor \frac{n}{3}\rfloor }\omega
_{i,1}(n|a)+I_{-2}(a)\Gamma (4/3)\sum_{k=1}^{\lfloor \frac{n+1}{3}\rfloor
}\omega _{i,2}(n|a) \\
&&+Ai(a)\sum_{\mu =0}^{2}\{\sum_{k=0}^{\infty }\Omega _{k,\mu }(n|a)\,\sigma
_{k,\mu }(a)+\sum_{k=1}^{\lfloor \frac{n+1-\mu }{3}\rfloor }\omega _{k,\mu
}(n|a)\,\widehat{\sigma }_{k,\mu }^{\,}(a)\} \\
&&-Ai^{\,\prime }(a)\sum_{\mu =0}^{2}\{\sum_{k=0}^{\infty }\Omega _{k,\mu
}(n|a)\,\,\sigma _{k,\mu }^{\,\prime }(a)+\sum_{k=1}^{\lfloor \frac{n+1-\mu 
}{3}\rfloor }\omega _{k,\mu }(n|a)\,\widehat{\sigma }_{k,\mu }^{\,\prime
}(a)\}.
\end{eqnarray*}

\bigskip

\begin{center}
\appendix Appendix V
\end{center}

Another method of computing the $\mathcal{I}_{n},i_{n},i_{n}^{\,\prime }$
integrals follows if we note that \ $\mathcal{I}_{n}(a)$ has the form 
\begin{equation*}
\mathcal{I}_{n}(a)=p_{n}(a)\,Ai(a)^{2}+q_{n}(a)\,Ai^{\,\,\prime
}(a)^{2}+r_{n}(a)\,Ai(a)\,Ai^{\,\,\prime }(a),
\end{equation*}%
where $p_{n}(a),$ $q_{n}(a),$ and $r_{n}(a)$ are polynomials. \ As a result
we can f\thinspace ind recurrence relations for these polynomials using
equation (26). \ We get%
\begin{eqnarray*}
2(2n-1)p_{n}(a)-n(n-1)(n-2)p_{n-3}(a) &=&-(n-1)a^{n}, \\
2(2n-1)q_{n}(a)-n(n-1)(n-2)q_{n-3}(a) &=&-na^{n-1}, \\
2(2n-1)r_{n}(a)-n(n-1)(n-2)r_{n-3}(a) &=&n(n-1)a^{n-2}.
\end{eqnarray*}%
It also follows that $i_{n}(a)$ and $i_{n}^{\,\,\prime }(a)$ have the same
form as $\mathcal{I}_{n}(a)$ and we have in summary

\begin{eqnarray*}
\mathcal{I}_{n}(a) &=&p_{n}(a)Ai(a)^{2}+q_{n}(a)Ai^{\,\,\prime
}(a)^{2}+r_{n}(a)Ai(a)Ai^{\,\,\prime }(a), \\
(n+1)i_{n}(a) &=&-[2p_{n+1}(a)+a^{n+1}]\,Ai(a)^{2} \\
&&-2q_{n+1}(a)\,Ai^{\,\,\prime }(a)^{2}-2r_{n+1}(a)\,Ai(a)Ai^{\,\,\prime
}(a), \\
(n+1)i_{n}^{\,\prime }(a)
&=&-2p_{n+2}(a)\,Ai(a)^{2}-[2q_{n+2}(a)+a^{n+1}]\,Ai^{\,\,\prime }(a)^{2} \\
&&-2r_{n+2}(a)\,Ai(a)Ai^{\,\,\prime }(a).
\end{eqnarray*}%
The f\thinspace irst few of these polynomials are given in the Table 7.

\begin{center}
\bigskip 
\begin{table}[H] \centering%
\caption{The polynomials p,q and r\label{eight}}%
\end{table}%

\begin{tabular}{|c|c|c|c|}
\hline
$n$ & $p_{n}(a)$ & $q_{n}(a)$ & $r_{n}(a)$ \\ \hline
$0$ & $-1/2$ & $0$ & $0$ \\ \hline
$1$ & $0$ & $-1/2$ & $0$ \\ \hline
$2$ & $-a^{2}/6$ & $-a/3$ & $1/3$ \\ \hline
$3$ & $-(3+2a^{3})/10$ & $-3a^{2}/10$ & $3a/5$ \\ \hline
$4$ & $-3a^{4}/14$ & $-2(3+a^{3})/7$ & $6a^{2}/7$ \\ \hline
$5$ & $-(5a^{2}+2a^{5})/9$ & $-(20a+5a^{4})/18$ & $10(1+a^{3})/9$ \\ \hline
\end{tabular}
\end{center}

It should be noted that the $p,q,r$ polynomials are perhaps more fundamental
than the integrals $i_{n}(a),$ $i_{n}^{\,\,\prime }(a),$ and $\mathcal{I}%
_{n}(a)$ themselves since they form their common basis.

\begin{center}
\appendix Appendix VI
\end{center}

Expressing the $i_{k}(a)$ and $i_{k}^{\,\prime }(a)$ integrals in terms of $%
\mathcal{I}_{k}(a)$ we have $(n\geq 1)$ where $\Xi ^{(n)}(-1,0),$ $\Lambda
^{(n)}(-1,0)$, and $\varrho ^{(n)}(-1,0)$ are hereafter denoted $\Xi ^{(n)},$
$\Lambda ^{(n)}$ and $\varrho ^{(n)}$ we have 
\begin{eqnarray*}
J_{n}(a) &=&\sum_{k=0}^{\infty }\{%
\begin{array}{c}
\,\varrho \mathfrak{\,}^{(k+n)}\mathcal{I}_{k}(a)-\tfrac{\Xi ^{(k+n)}}{(k+1)}%
\,[2\mathcal{I}_{k+1}(a)+a^{k+1}Ai(a)^{2}] \\ 
-\tfrac{\Lambda ^{(k+n)}}{(k+1)}[2\mathcal{I}_{k+2}(a)+a^{k+1}Ai^{\,\prime
}(a)^{2}]%
\end{array}%
\}/(k+n)! \\
&&+\left\{ \varrho \mathfrak{\,}^{(n-1)}\,\mathcal{I}_{-1}(a)+\,\Xi
^{(n-1)}i_{-1}(a)+\Lambda ^{(n-1)}i_{-1}^{\,\prime }(a)\right\} /(n-1)! \\
&&+\sum_{k=2}^{n}\{%
\begin{array}{c}
\,\varrho \mathfrak{\,}^{(n-k)}\mathcal{I}_{-k}(a) \\ 
+\tfrac{\Xi ^{(n-k)}}{(k-1)}\,[2\mathcal{I}_{-k+1}(a)+Ai(a)^{2}/a^{k-1}] \\ 
+\tfrac{\Lambda ^{(n-k)}}{(k-1)}[2\mathcal{I}_{-k+2}(a)+a^{-k+1}Ai^{\,\prime
}(a)^{2}]%
\end{array}%
\}/(n-k)!
\end{eqnarray*}%
Rearranging this expression and gathering common $\mathcal{I}_{k}(a)$ terms
we have 
\begin{eqnarray*}
J_{n}(a) &=&-Ai(a)^{2}\sum_{k\,=\,-n}^{\infty }\tfrac{\Xi ^{(k+n)}a^{k+1}}{%
(k+1)(n+k)!}\mathbf{\Delta }_{k,-1}-Ai^{\,\prime }(a)^{2}\sum_{k=-n}^{\infty
}\tfrac{\Lambda ^{(k+n)}a^{k+1}}{(k+1)(n+k)!}\mathbf{\Delta }_{k,-1} \\
&&+\left\{ \varrho \mathfrak{\,}^{(n)}/n!\right\} \,\mathcal{I}%
_{0}(a)+\left\{ \varrho \mathfrak{\,}^{(n+1)}/(n+1)!-2\,\Xi
^{(n)}/n!\right\} \,\mathcal{I}_{1}(a) \\
&&+\sum_{k=2}^{\infty }\{\,\frac{\varrho \mathfrak{\,}^{(k+n)}}{(k+n)!}-2%
\frac{\Xi ^{(k+n-1)}}{k(k+n-1)!}-2\frac{\Lambda ^{(k+n-2)}}{(k-1)(k+n-2)!}\}%
\mathcal{I}_{k}(a) \\
&& \\
&&+2\left\{ \Lambda ^{(n-2)}/(n-2)!\right\} \,\mathcal{I}_{0}(a)+\left\{
\,\Xi ^{(n-1)}/(n-1)!\right\} i_{-1}(a) \\
&&+\left\{ \Lambda ^{(n-1)}/(n-1)!\right\} i_{-1}^{\,\prime }(a) \\
&&+\sum_{k=1}^{n-2}\,\{\frac{\varrho \mathfrak{\,}^{(n-k)}}{(n-k)!}+2\frac{%
\Xi ^{(n-k-1)}}{k(n-k-1)!}+2\frac{\Lambda ^{(n-k-2)}}{(k+1)(n-k-2)!}\}%
\mathcal{I}_{-k}(a) \\
&&+\left\{ \varrho \mathfrak{\,}^{(0)}\right\} \,\mathcal{I}%
_{-n}(a)+\{\varrho \mathfrak{\,}^{(1)}+2\,\mathbf{\Delta }_{n,1}\,\frac{\Xi
^{(0)}}{(n-1)}\,\}\mathcal{I}_{-n\,+1}(a),
\end{eqnarray*}%
with $\mathbf{\Delta }_{k,\,\,j}\equiv 1-\delta _{k,\,\,j\text{ }}$where $%
\delta _{k,\,\,j}$ is the Kronecker delta function. \ Since the integrals $%
\mathcal{I}_{k}(a)$ and $\mathcal{I}_{-k}(a)$ are related to the Airy
functions by (27) and (30) we have for $J_{n}(a)$ in the cases $(n\geq 4)$ 
\begin{eqnarray*}
J_{n}(a) &=&Ai(a)^{2}\{%
\begin{array}{c}
-\,\Xi _{n}(a)+\sum_{\mu =0}^{2}\sum_{k=0}^{\infty }\Psi _{k,\mu
}(n|a)\,B_{k,\mu }^{(0)} \\ 
+\sum_{\mu =0}^{2}\sum_{k=0}^{\lfloor \frac{n-2-\mu }{3}\rfloor }\psi
_{k,\mu }(n|a)\,\beta _{k,\mu }^{(0)}%
\end{array}%
\} \\
&&+Ai^{\,\prime }(a)^{2}\{%
\begin{array}{c}
-\,\Lambda _{n}(a)+\sum_{\mu =0}^{2}\sum_{k=0}^{\infty }\Psi _{k,\mu
}(n|a)\,B_{k,\mu }^{(1)} \\ 
+\sum_{\mu =0}^{2}\sum_{k=0}^{\lfloor \frac{n-2-\mu }{3}\rfloor }\psi
_{k,\mu }(n|a)\,\beta _{k,\mu }^{(1)}%
\end{array}%
\} \\
&&+Ai(a)Ai^{\,\prime }(a)\{%
\begin{array}{c}
\sum_{\mu =0}^{2}\sum_{k=0}^{\infty }\Psi _{k,\mu }(n|a)\,B_{k,\mu }^{(2)}
\\ 
+\sum_{\mu =0}^{2}\sum_{k=0}^{\lfloor \frac{n-2-\mu }{3}\rfloor }\psi
_{k,\mu }(n|a)\,\beta _{k,\mu }^{(2)}%
\end{array}%
\} \\
&&+\mathcal{I}_{1}\{\tfrac{\rho ^{(n+1)}}{(n+1)!}-2\tfrac{\Xi ^{(n)}}{(n)!}%
\}+\mathcal{I}_{0}\{\tfrac{\rho ^{(n)}}{n!}+2\tfrac{\Lambda ^{(n-2)}}{(n-2)!}%
\} \\
&&+\tfrac{12}{\sqrt{\pi }}\mathcal{I}_{-1}\{\sum_{k=0}^{\lfloor \frac{n-3}{3}%
\rfloor }\psi _{k,1}(n|a)\} \\
&&+\tfrac{12}{\Gamma (5/6)}\mathcal{I}_{-2}\{\sum_{k=0}^{\lfloor \frac{n-4}{3%
}\rfloor }\psi _{k,2}(n|a)\}+\tfrac{12}{\Gamma (1/6)}\mathcal{I}%
_{-3}\{\sum_{k=0}^{\lfloor \frac{n-2}{3}\rfloor }\psi _{k,3}(n|a)\} \\
&&+\tfrac{\Xi ^{(n-1)}}{(n-1)!}i_{-1}+\tfrac{\Lambda ^{(n-1)}}{(n-1)!}%
i_{-1}^{\,\prime }+\varrho ^{(0)}\mathcal{I}_{-n}+\{\varrho ^{(1)}+2\Delta
_{n,1}\tfrac{\Xi ^{(0)}}{(n-1)}\}\mathcal{I}_{-n+1}
\end{eqnarray*}%
where%
\begin{eqnarray*}
\Xi _{n}(a) &=&\sum_{k\,=\,-n}^{\infty }\tfrac{\Xi ^{(k+n)}a^{k+1}}{%
(k+1)(n+k)!}\mathbf{\Delta }_{k,-1}, \\
\Lambda _{n}(a) &=&\sum_{k=-n}^{\infty }\tfrac{\Lambda ^{(k+n)}a^{k+1}}{%
(k+1)(n+k)!}\mathbf{\Delta }_{k,-1},
\end{eqnarray*}%
and the terms $\Psi _{k,\mu }$ and $\psi _{k,\mu }$ are given by%
\begin{eqnarray*}
\Psi _{k,\mu }(n|a) &=&\tfrac{\Gamma (3k+\mu +1)}{12^{k+1}\Gamma (k+\mu
/3+5/6)} \\
&&\cdot \{\tfrac{\varrho ^{(n+3k+\mu )}}{(n+3k+\mu )!}-2\tfrac{\Xi
^{(n+3k+\mu -1)}}{^{(3k+\mu )\,(n+3k+\mu -1)!}}-2\tfrac{\Lambda ^{(n+3k+\mu
-2)}}{^{(3k+\mu -1)\,(n+3k+\mu -2)!}}\}, \\
\psi _{k,\mu }(n|a) &=&\tfrac{12^{k-1}\Gamma (k+\mu /3+1/6)}{\Gamma (3k+\mu )%
} \\
&&\cdot \{\tfrac{\varrho ^{(n-3k-\mu )}}{(n-3k-\mu )!}+2\tfrac{\Xi
^{(n-3k-\mu -1)}}{^{(3k+\mu )\,(n-3k-\mu -1)!}}+2\tfrac{\Lambda ^{(n-3k-\mu
-2)}}{^{(3k+\mu +1)\,(n-3k-\mu -2)!}}\}.
\end{eqnarray*}


\begin{thebibliography}{99}
\bibitem{airynbs} M. Abramowitz and I. A. Stegun, editors. \ \textit{%
Handbook of Mathematical Functions With Formulas, Graphs, And Mathematical
Tables}, New York: Dover, p. 446 (1992).

\bibitem{mathphy} B. J. Laurenzi, \textit{An Analytic Solution to the
Thomas-Fermi Equations,} J. Math. Phys. \textbf{31}, p. 2536 (1990).

\bibitem{laur} B. J. Laurenzi, A preliminary account of the present work is
found in \textit{Logarithmic Integrals of Airy Functions},
arXiv:1211.0705v1.pdf, .

\bibitem{vallee} O. Vallee and M. Soares, \textit{Airy Functions and
Applications to Physics, }second edition\textit{, }Imperial College Press,
p. 19, (2010).

\bibitem{nbs} Ref. 1, p. 478. Table 10.13.Ref. 8, http://dlmf.nist.gov/11.10

\bibitem{Olver} F. W. J. Olver, \textit{Airy and Related Functions},
[DLMF]NIST Digital Library of Mathematical Functions.
http://dlmf.nist.gov/9.9E8

\bibitem{laurpq} B. J. Laurenzi, \textit{Polynomials Associated with the
Higher Derivatives of the Airy Functions }$Ai(z)$ \textit{and }$Ai^{\,\prime
}(z)$, arXiv:1110.2025.

\bibitem{nbsscorer} NIST Digital Library of Mathematical Functions, \textit{%
Airy and Related Functions}, http://dlmf.nist.gov/9.12.

\bibitem{Aspnes} D. E. Aspnes, \textit{Electric-f\thinspace ields
Ef\thinspace f\thinspace ects on Optical Absorption Near Thresholds in
Solids,} Phys. Rev. \textbf{147}, (1966) pp. 554-566.

\bibitem{Reid} W. H. Reid, \textit{Integral Representation for Products of
Airy Functions,} Z. Angew. Math. Phys. \textbf{46, }(1995) pp. 159-170.
\end{thebibliography}
\end{document}